\documentclass{amsart}
\usepackage{amsmath,amssymb,amsthm,mathtools}
\usepackage[margin=1in]{geometry}
\usepackage[hidelinks]{hyperref}
\hypersetup{
 pdftitle={schubmult: Symbolic algorithms and input-sensitive bounds for Schubert multiplication},
 pdfauthor={Matthew J. Samuel}
}
\makeatletter
\renewcommand{\paragraph}{\@startsection{paragraph}{4}%
  \z@{.5\linespacing\@plus.3\linespacing}{.25\linespacing}%
  {\normalfont\bfseries}}
\makeatother
\newtheorem{theorem}{Theorem}[section]
\newtheorem{proposition}[theorem]{Proposition}
\newtheorem{lemma}[theorem]{Lemma}
\newtheorem{corollary}[theorem]{Corollary}
\theoremstyle{definition}

\newtheorem{method}[theorem]{Method}

\theoremstyle{remark}

\numberwithin{equation}{section}
\newcommand{\Sch}{\mathfrak S}
\newcommand{\code}{\mathfrak c}
\newcommand{\id}{\mathrm{id}}
\newcommand{\Des}{\mathrm{Des}}
\newcommand{\fall}[2]{(#1)_{\underline{#2}}}
\newcommand{\alg}{\texttt{schubmult}}

\title[\alg: Schubert multiplication]{\alg: Symbolic algorithms and input-sensitive bounds for Schubert multiplication}
\author{Matthew J. Samuel}
\email{matthematics@gmail.com}
\date{}

\begin{document}
\begin{abstract}
We present \alg, a GPL-licensed Python system for exact multiplication of
ordinary, double, quantum, and quantum double Schubert polynomials, as
well as independent-alphabet triple and quantum triple products. Its
central algorithm expands one factor into filtered column paths and
couples them to Pieri computations on the other factor. Dynamic
programming shares intermediate states, and the algorithm returns every
Schubert-basis coefficient as an exact arithmetic circuit. We prove
input-sensitive operation bounds for the entire product, not just a
selected coefficient: for inputs of minimal rank $n$, the classical
algorithms have complexity $2^{O(n\log n)}$ and polynomial complexity in
$n$ for fixed Lehmer-code lengths. The quantum algorithms have the same
general bound and polynomial complexity when, in addition, the number
of nonzero entries in the Lehmer code of $v^{-1}$ is fixed. We describe the
implementation and distinguish it from the algorithms analyzed, verify
computations against independent routines and exact identities, and
report reproducible timings for flag and Grassmannian workloads.
Parabolic quantum products are treated separately via the
Peterson--Woodward comparison.
\end{abstract}
\maketitle

\section{Introduction}

\alg\ \cite{schubmult} computes ordinary, double, quantum, and quantum
double Schubert products, as well as independent-alphabet triple and
quantum triple products. It also provides parabolic quantum and
parabolic quantum double multiplication. This paper presents the
multiplication architecture and analyzes its complexity in terms of
the input permutations.

For permutations $u,v$ with finite support, consider the expansion
\begin{equation}\label{eq:product}
 \Sch_u(x;y)\Sch_v(x;z)
   =\sum_w c_{u,v}^w(y;z)\Sch_w(x;y).
\end{equation}
Here $x,y,z$ are independent infinite alphabets, and the coefficients belong
to $\mathbb Z[y,z]$.
We call \eqref{eq:product} the \emph{triple} product.
At $z=y$ its coefficients are the equivariant Schubert structure
constants of complete flag varieties.
At $y=z=0$ they are the ordinary Schubert structure constants, which
are nonnegative integers.
For $y\ne z$ the coefficients are more general than both.
The basic multiplication task of \alg\ is to compute all the coefficients in
\eqref{eq:product}, not just a coefficient specified by a third permutation.
We study full-product computation as an exact symbolic-computation
problem: the output is a Schubert-basis dictionary whose coefficients are
arithmetic circuits, not necessarily expanded polynomials. This output
model keeps shared subexpressions and avoids requiring symbolic
simplification or zero testing.

\alg\ computes such products by matching two
computations. An elementary-symmetric expansion of the second factor
provides a graph of signed labels. Pieri's rule multiplies the first factor
by each elementary factor. Dynamic programming combines these operations
without expanding the final coefficients into monomials.
The expansion formulas and correctness of this architecture are established
in \cite[Section~4 and Section~7]{samuelquantum}.
The contribution here is an input-sensitive analysis of the full-product
algorithms: a justified path filter and explicit counts for every
intermediate state and transition, including the additional quantum
states. The operation bounds are separated from the measured costs of
the released symbolic implementation.
We analyze this architecture through bounds in parameters recoverable
from the input codes, rather than by enumerating the ambient symmetric
group.
We use ``classical'' only to mean ``non-quantum''; it includes the
double and triple products.
In product names, ``quantum'' denotes the non-equivariant product,
``quantum double'' retains one shared equivariant alphabet, and
``quantum triple'' retains two independent coefficient alphabets.
These are distinct multiplication tasks, even though the same architecture
and state-counting argument handle all three.

We use the distinct-position column paths of
\cite[Definition~4.4 and Theorem~4.5]{samuelquantum};
their restriction follows from the cancellation theorem proved there.
The arithmetic-circuit model and the choice of elementary-weight evaluator
are distinguished from literal symbolic implementation costs in
Section~\ref{sec:implementation}.

\paragraph{Software availability}
The package provides Python ring objects and command-line programs,
with symbolic arithmetic based on SymEngine and SymPy; where a compatible
wheel is available, it also uses a compiled C++ backend.
It is distributed under the GPL-3.0-or-later license, and
\texttt{pip install schubmult} installs the latest release.
Each release declares its supported Python versions in its package
metadata; these may change between releases.
Source code is available at
\url{https://github.com/matthematics/schubmult},
with documentation and usage examples at
\url{https://matthematics.github.io/schubmult/}.
The benchmarked package is the non-editable PyPI release
\texttt{schubmult==5.3.0}, installed in a CPython~3.14.10 virtual
environment; the release source is available at
\url{https://github.com/matthematics/schubmult/tree/5.3.0}.
That release's metadata requires Python~3.10 or later; this is a
release-specific minimum, not a promise of support for every future
Python version.
The benchmark script and raw measurements accompany the manuscript; the
experimental protocol and environment are detailed in
Section~\ref{sec:practical}.
Appendix~\ref{app:grothendieck} describes the Grothendieck interfaces,
including the conjectural equivariant quantum extension.

\paragraph{Buch's \texttt{schubmult} and \texttt{lrcalc}}
Buch's Littlewood--Richardson Calculator (\texttt{lrcalc}) \cite{lrcalc}
is integrated into SageMath \cite{sagemath}. Its main purpose is the
computation of ordinary Grassmannian Littlewood--Richardson coefficients
and products of Schur functions; Sage's \texttt{lrcalc.mult} interface
also computes products in the quantum cohomology of a Grassmannian
when its \texttt{quantum}, \texttt{maxrows}, and \texttt{level}
arguments are supplied \cite{sagelrcalcdocs}.
Its specialized routines are a natural alternative for these tasks.
The measurements in Section~\ref{sec:practical} concern \alg\ and
do not establish a performance ranking among packages.

The \texttt{lrcalc} distribution also includes Buch's earlier C program
\texttt{schubmult}, which computes products of ordinary Schubert
polynomials and was for many years the standard tool for doing so.
It expands one factor into monomials by the transition recursion of
Lascoux and Sch\"utzenberger \cite{lasschubschur} and multiplies the
monomials into the other factor one variable at a time by Monk's rule
\cite{monk}.
The present Python package is independent and does not include Buch's
implementation; its main kernels instead match an elementary-symmetric
expansion of one factor against Pieri multiplication of the other, as
described below.
For comparison, our package contains a separately written implementation
of the same transition--Monk algorithm, together with a heuristic
choosing between the two methods for ordinary products.
These routines are available as Python functions, not through the
command-line interface, and are not analyzed here.

\paragraph{Organization}
Section~\ref{sec:completion} constructs the dominant completion and
proves the label filter. Sections~\ref{sec:columns}--\ref{sec:pieri}
derive the filtered elementary expansion and the multiplication
algorithm. Sections~\ref{sec:states}--\ref{sec:complexity} establish
the classical bounds, and Section~\ref{sec:quantum} treats quantum,
quantum double, and quantum triple multiplication.
Section~\ref{sec:dictionary} gives the cohomological interpretation of
ordinary, double, quantum, and quantum double products, distinguishing
complete flags from partial flags.
Section~\ref{sec:implementation} discusses the
implementation, computational checks, and parabolic output.
Section~\ref{sec:practical} reports reproducible timings on
application-based workloads and a controlled family of flag products.
The appendix documents the Grothendieck interface.

\subsection{Cohomological dictionary}\label{sec:dictionary}

\paragraph{Complete flags}
The ordinary, double, quantum, and quantum double products computed here
have the following interpretations in the cohomology of a complete flag
variety. The Schubert polynomial $\Sch_w(x;0)$ represents the ordinary
Schubert class $\sigma_w$; its double version $\Sch_w(x;y)$ represents
the $T$-equivariant class $\sigma_w^T$, where $y$ records the characters
of the diagonal torus. The quantum Schubert polynomial
$\widetilde\Sch_w(x;0)$ represents the non-equivariant quantum Schubert
class in $QH^*(\mathrm{Fl}_N)$. The quantum double Schubert polynomial
$\widetilde\Sch_w(x;y)$ represents the equivariant quantum Schubert
class in $QH_T^*(\mathrm{Fl}_N)$
\cite{billey1999kostant,fomin1997quantum,lam2014quantum}.
In each product below, both factors use the same coefficient alphabet:
\[
\begin{aligned}
 \Sch_u(x;0)\Sch_v(x;0)
   &=\sum_w c_{u,v}^w(0;0)\Sch_w(x;0),
       && H^*(\mathrm{Fl}_N),\\
 \Sch_u(x;y)\Sch_v(x;y)
   &=\sum_w c_{u,v}^w(y;y)\Sch_w(x;y),
       && H_T^*(\mathrm{Fl}_N),\\
 \widetilde\Sch_u(x;0)\widetilde\Sch_v(x;0)
   &=\sum_{w,\delta}q^\delta c_{u,v}^{w,\delta}(0;0)
                    \widetilde\Sch_w(x;0),
       && QH^*(\mathrm{Fl}_N),\\
 \widetilde\Sch_u(x;y)\widetilde\Sch_v(x;y)
   &=\sum_{w,\delta}q^\delta c_{u,v}^{w,\delta}(y;y)
                    \widetilde\Sch_w(x;y),
       && QH_T^*(\mathrm{Fl}_N).
\end{aligned}
\]
Thus ordinary and double products yield ordinary and equivariant Schubert
structure constants, respectively. Quantum products yield genus-zero
three-point Gromov--Witten invariants; quantum double products yield
equivariant Gromov--Witten polynomials.
Here $q^\delta=\prod_{i=1}^{N-1}q_i^{\delta_i}$
records the curve degree in the simple-coroot basis. These are four
complete-flag cases of the multiplication architecture.
Parabolic quantum and parabolic quantum double products require a
separate interpretation and comparison, described below.

For clarity, the finite-rank presentations use $x_1,\ldots,x_N$ and,
equivariantly, $y_1,\ldots,y_N$. With $e_j$ the elementary symmetric
polynomial, the classical and equivariant rings are
\[
 H^*(\mathrm{Fl}_N)
   \cong \mathbb Z[x_1,\ldots,x_N]/(e_j(x))_{j=1}^N,\qquad
 H_T^*(\mathrm{Fl}_N)
   \cong \mathbb Z[y_1,\ldots,y_N,x_1,\ldots,x_N]/
        (e_j(x)-e_j(y))_{j=1}^N.
\]
For the quantum presentations, let $C_N(x;q)$ be the tridiagonal matrix
with diagonal entries $x_i$, superdiagonal entries $-1$, and
subdiagonal entries $q_i$, and define $\mathcal E_j^{(N)}(x;q)$ by
\[
 \det(C_N(x;q)-tI)
   =\sum_{j=0}^N(-t)^{N-j}\mathcal E_j^{(N)}(x;q).
\]
Then
\[
 \begin{aligned}
 QH^*(\mathrm{Fl}_N)
   &\cong \mathbb Z[x_1,\ldots,x_N,q_1,\ldots,q_{N-1}]
       /(\mathcal E_j^{(N)}(x;q))_{j=1}^N,\\
 QH_T^*(\mathrm{Fl}_N)
   &\cong \mathbb Z[y_1,\ldots,y_N,x_1,\ldots,x_N,
                    q_1,\ldots,q_{N-1}]
       /(\mathcal E_j^{(N)}(x;q)-e_j(y))_{j=1}^N.
 \end{aligned}
\]
These presentations also specify how to interpret stable polynomial
calculations at finite rank: map to the rank-$N$ quotient, whose Schubert
basis is indexed by $S_N$. In the complete-flag polynomial grading,
$\deg x_i=\deg y_i=1$ and $\deg q_i=2$; these are half the usual
cohomological degrees \cite{lam2014quantum}.

\paragraph{Partial flags and classical pullback}
Fix positive block sizes $\mathbf b=(b_1,\ldots,b_{s+1})$ summing to $N$,
and set $a_0=0$ and $a_j=b_1+\cdots+b_j$. The corresponding partial flag
variety is
\[
 F_{\mathbf b}=\mathrm{Fl}(a_1,\ldots,a_s;N)=GL_N/P.
\]
Its Schubert classes are indexed by the minimal right-coset
representatives $W^P$ for
$S_N/W_P$, where $W_P=S_{b_1}\times\cdots\times S_{b_{s+1}}$.
Equivalently, $w\in W^P$ is increasing within each block of positions.
The projection $\pi:\mathrm{Fl}_N\to F_{\mathbf b}$ induces injective
ring maps on ordinary and equivariant cohomology, with
$\pi^*(\sigma_w^P)=\sigma_w$ and
$\pi^*(\sigma_w^{P,T})=\sigma_w^T$ for $w\in W^P$.
Thus the same ordinary and double Schubert polynomials represent
the pulled-back classical classes \cite[Section~6]{lam2014quantum}.

\paragraph{Parabolic quantum and parabolic quantum double products}
The analogous identification of Schubert bases does \emph{not} in general
extend $\pi^*$ to a quantum ring homomorphism. The rings
$QH^*(F_{\mathbf b})$ and $QH_T^*(F_{\mathbf b})$ carry the parabolic
quantum and parabolic quantum double products, respectively.
Their polynomial representatives are the
\emph{parabolic} quantum Schubert polynomials
$\widetilde\Sch_w^P(x;0)$ and quantum double Schubert polynomials
$\widetilde\Sch_w^P(x;y)$. These depend on $P$ and are generally
different from the complete-flag polynomials
$\widetilde\Sch_w(x;0)$ and $\widetilde\Sch_w(x;y)$ with the same index
\cite{ciocan1999partial,lam2014quantum}.

Write $Q_1,\ldots,Q_s$ for the parabolic quantum parameters, to distinguish
them from the complete-flag parameters $q_1,\ldots,q_{N-1}$.
They record the degrees of the Schubert curves corresponding to the
block boundaries $a_1,\ldots,a_s$, a basis of $H_2(F_{\mathbf b},\mathbb Z)$.
As modules, the two rings are
\[
 \begin{aligned}
 QH^*(F_{\mathbf b})
   &\cong\bigoplus_{w\in W^P}\mathbb Z[Q_1,\ldots,Q_s]\,\sigma_w^P,\\
 QH_T^*(F_{\mathbf b})
   &\cong\bigoplus_{w\in W^P}\mathbb Z[y_1,\ldots,y_N,Q_1,\ldots,Q_s]\,
                    \sigma_w^{P,T}.
 \end{aligned}
\]
The parabolic quantum double product is
\[
 \sigma_u^{P,T}\star_P\sigma_v^{P,T}
   =\sum_{\substack{w\in W^P\\\delta\in\mathbb Z_{\geq0}^{\,s}}}
       Q^\delta C_{u,v}^{w,\delta;P}(y)\,\sigma_w^{P,T},
 \qquad Q^\delta=\prod_{j=1}^s Q_j^{\delta_j}.
\]
Here $C_{u,v}^{w,\delta;P}(y)$ is the equivariant three-point
Gromov--Witten invariant of $F_{\mathbf b}$; setting $y=0$ gives the
parabolic quantum product, and setting $Q=0$ gives the classical
equivariant product. In the half-cohomological grading,
\[
 \deg Q_j=a_{j+1}-a_{j-1}=b_j+b_{j+1},
 \qquad
 \deg_y C_{u,v}^{w,\delta;P}(y)
   =\ell(u)+\ell(v)-\ell(w)
      -\sum_{j=1}^s(b_j+b_{j+1})\delta_j
\]
for a nonzero coefficient \cite[Section~6.2]{lam2014quantum}.
In particular, the quantum parameter of $\operatorname{Gr}(k,N)$
has polynomial degree $N$, not $2$. When all blocks have size $1$,
the complete-flag dictionary is recovered.

The interfaces \texttt{QPSx} and \texttt{QPDSx}, supplied with these
block sizes, compute parabolic quantum and parabolic quantum double
products, respectively; the latter uses the same $y$-alphabet in both
factors. Their printed \texttt{q\_j} denotes $Q_j$.
They obtain parabolic coefficients from complete-flag coefficients by
the Peterson--Woodward comparison, which uses a distinguished lift of
the curve degree and a modified output permutation
\cite{woodward2005d}, \cite[Proposition~2.10]{huang2015equivariant}.
This is not merely restriction to the labels $W^P$, nor is it obtained
by simply setting unwanted full-flag quantum parameters to zero.
Section~\ref{sec:parabolic-output} describes the comparison implemented
by the package. The subsequent \texttt{kill\_ideal()} operation imposes
the finite-rank quotient on the already parabolic output; it does not
perform the passage from a complete-flag product to the corresponding
parabolic product.

For example, take blocks $(2,2)$, so $F_{\mathbf b}=\operatorname{Gr}(2,4)$.
The products with the same input labels in the two finite quantum rings
are
\[
 \begin{aligned}
 \sigma_{1324}\star\sigma_{3412}
   &=q_2\sigma_{3142},
       &&QH^*(\mathrm{Fl}_4),\\
 \sigma_{1324}^P\star_P\sigma_{3412}^P
   &=Q_1\sigma_{1324}^P,
       &&QH^*(\operatorname{Gr}(2,4)).
 \end{aligned}
\]
The index $3142$ is not in $W^P$: discarding nonminimal indices would
erase the nonzero parabolic answer. The comparison changes both the
index and the interpretation of the degree, with $\deg q_2=2$ but
$\deg Q_1=4$. Section~\ref{sec:examples} gives this parabolic computation
and its parabolic quantum double version in the Python interface.

\paragraph{Scope}
The geometric dictionary uses a single equivariant alphabet shared by
the two factors. The product with independent alphabets $y$ and $z$
from \eqref{eq:product}, and its quantum analogue, is a formal
coefficient-refined triple product; when $z\ne y$ it is not
multiplication in one ordinary or equivariant cohomology ring.

\subsection{Examples of use}\label{sec:examples}

\paragraph{Python ring interface}
The ring objects accept permutations in one-line notation.
Multiplication returns a Schubert-basis element, whose coefficients
can be accessed using \texttt{Permutation} keys; \texttt{expand()}
converts it to an explicit polynomial. Trailing fixed points are
omitted in displayed permutation labels.
\begin{verbatim}
>>> from schubmult import Sx, Permutation
>>> p = Sx([3, 1, 2]) * Sx([1, 3, 2])
>>> p
Sx((3, 2, 1)) + Sx((4, 1, 2, 3))
>>> p[Permutation([3, 2, 1])]
1
>>> p.expand()
x_1**2*x_2 + x_1**3
\end{verbatim}
The double ring \texttt{DSx} also supports independent coefficient
alphabets. Its second argument selects the coefficient alphabet of
that factor; omitting it uses $y$.
\begin{verbatim}
>>> from schubmult import DSx
>>> from schubmult.abc import z
>>> DSx([1, 3, 2]) * DSx([2, 1, 3], z)
(y_1 - z_1)*DSx((1, 3, 2), y) + DSx((2, 3, 1), y) + DSx((3, 1, 2), y)
\end{verbatim}
The quantum ring \texttt{QSx} and quantum double ring \texttt{QDSx}
compute non-equivariant and equivariant quantum products, respectively:
\begin{verbatim}
>>> from schubmult import QSx, QDSx
>>> QSx([2, 1])**2
q_1 + QSx((3, 1, 2))
>>> QDSx([2, 1])**2
q_1 + (-y_1 + y_2)*QDSx((2, 1), y) + QDSx((3, 1, 2), y)
\end{verbatim}
The parabolic quantum ring \texttt{QPSx} and parabolic quantum double
ring \texttt{QPDSx} take the block sizes. The following computes
$\sigma_{1324}\star\sigma_{3412}=q_1\sigma_{1324}$ in
$QH^*(\operatorname{Gr}(2,4))$ and its equivariant counterpart in
$QH_T^*(\operatorname{Gr}(2,4))$; \texttt{kill\_ideal()} imposes the
finite-rank quotient after parabolic multiplication.
\begin{verbatim}
>>> from schubmult import QPSx, QPDSx
>>> A = QPSx(2, 2)
>>> (A([1, 3, 2, 4]) * A([3, 4, 1, 2])).kill_ideal()
q_1*QPSx(2, 2)((1, 3, 2))
>>> B = QPDSx(2, 2)
>>> (B([1, 3, 2, 4]) * B([3, 4, 1, 2])).kill_ideal()
q_1*QPDSx(2, 2)((1, 3, 2), y) + (-y_1 - y_2 + y_3 + y_4)*QPDSx(2, 2)((3, 4, 1, 2), y)
\end{verbatim}
The parabolic quantum double ring also supports parabolic quantum triple
products with independent coefficient alphabets. Passing $z$ to the second
factor gives an expansion in the $y$-basis:
\begin{verbatim}
>>> (B([3, 4, 1, 2]) * B([1, 3, 2, 4], z)).kill_ideal()
q_1*QPDSx(2, 2)((1, 3, 2), y) + (y_3 + y_4 - z_1 - z_2)*QPDSx(2, 2)((3, 4, 1, 2), y)
\end{verbatim}

\paragraph{Command-line tools}
Installing the Python package provides the commands below.
Permutations are space-separated integers, with a lone \texttt{-}
separating factors. The ordinary program prints the coefficient
followed by the permutation. With \texttt{--code}, both input and
output use Lehmer codes:
\begin{verbatim}
$ schubmult_py 3 1 2 - 1 3 2
1  (3, 2, 1)
1  (4, 1, 2, 3)
$ schubmult_py --code 2 - 0 1
1  [2, 1]
1  [3]
\end{verbatim}
The double, quantum, and quantum double programs print the permutation
first. For the double and quantum double programs, \texttt{--mixed-var}
selects independent $y,z$ alphabets; otherwise both factors use $y$.
\begin{verbatim}
$ schubmult_double 1 3 2 - 2 1 3 --mixed-var
(1, 3, 2)  y_1 - z_1
(2, 3, 1)  1
(3, 1, 2)  1
$ schubmult_q_double 2 1 - 2 1
()  q_1
(2, 1)  -y_1 + y_2
(3, 1, 2)  1
\end{verbatim}
Here \texttt{()} denotes the identity permutation. The arguments to
\texttt{--parabolic} are block sizes, so the parabolic examples above
can also be evaluated by
\begin{verbatim}
$ schubmult_q 1 3 2 4 - 3 4 1 2 --parabolic 2 2
(1, 3, 2)  q_1
$ schubmult_q_double 1 3 2 4 - 3 4 1 2 --parabolic 2 2
(1, 3, 2)  q_1
(3, 4, 1, 2)  -y_1 - y_2 + y_3 + y_4
\end{verbatim}
Combining \texttt{--parabolic} with \texttt{--mixed-var} computes
the parabolic quantum triple product:
\begin{verbatim}
$ schubmult_q_double 3 4 1 2 - 1 3 2 4 --parabolic 2 2 --mixed-var
(1, 3, 2)  q_1
(3, 4, 1, 2)  y_3 + y_4 - z_1 - z_2
\end{verbatim}

\subsection{Parameters and main bound}

Permutations are composed as functions, with $(ab)(i)=a(b(i))$.
Write
\[
 \code_i(w)=\#\{j>i:w(j)<w(i)\},\qquad
 \ell(w)=\sum_i\code_i(w),\qquad
 r_w=\max\{i:\code_i(w)>0\},
\]
where $r_{\id}=0$. Thus $r_w$ is the length of the trimmed Lehmer code
and the last descent position of $w$.
We take $n\geq2$ to be the least integer for which $u,v\in S_n$,
apart from this harmless lower bound of $2$.
It can be recovered from the input codes:
\[
 n=\max\bigl(\{2\}\cup\{i+\code_i(u):\code_i(u)>0\}
                     \cup\{i+\code_i(v):\code_i(v)>0\}\bigr).
\]
Put $r=\max(r_u,r_v)$, $h=r_v$, and $d=\ell(u)+\ell(v)$.
The case $v=\id$ requires no multiplication and will be set aside.

Section~\ref{sec:completion} constructs a dominant completion with column
heights $\theta_1\geq\cdots\geq\theta_L>0$.
It satisfies $\theta_1=h$ and $L\leq n-1$.
Define the following quantities, all determined or bounded by the input:
\begin{equation}\label{eq:parameters}
 \begin{split}
 M&=n+L,\\
 U&=\min\left\{\fall{M}{r},\binom{d+r}{r}\right\},\\
 V&=\fall{n}{h},\qquad
 Q_h(n)=\sum_{b=0}^{h}\fall{n-1}{b}.
 \end{split}
\end{equation}
Here $\fall{a}{b}=a!/(a-b)!$ and $\fall{a}{0}=1$.

\begin{theorem}[Input-parameter bound]\label{thm:main}
For $v\ne\id$, Method~\ref{method:product} computes the entire triple
product \eqref{eq:product}, with coefficient arithmetic circuits, in
\begin{equation}\label{eq:main}
 O\left(n^3+L M^3 U^2 V Q_h(n)\right)
\end{equation}
elementary and coefficient-ring operations.
The same bound, with the same proof, holds for the double ($z=y$) and
ordinary ($y=z=0$) products; see Proposition~\ref{prop:uniform}.
In particular, the bound is
\begin{equation}\label{eq:coarse}
 O\left(2^{2r}n^{2r+2h+4}\right)
   =2^{O(n\log n)}.
\end{equation}
For fixed $r_u,r_v$ it is polynomial in $n$.
If $v$ is dominant, the bound improves to
$O(n^3+L M^3 U^2)$.
\end{theorem}

The proof is in Section~\ref{sec:complexity}.
The exponent in \eqref{eq:coarse} is at most $4n\log_2 n+O(n)$.

\subsection{Quantum, quantum double, and quantum triple multiplication}

The same architecture computes quantum Schubert products, quantum double
Schubert products, and the quantum triple product \eqref{eq:qproduct}.
The quantum double case ($z=y$) gives the equivariant quantum
Littlewood--Richardson coefficients of complete flag varieties;
the quantum case ($y=z=0$) gives their non-equivariant counterparts.
Two further input parameters enter.
Let $L$ be the number of nonzero entries of $\code(v^{-1})$; it is
the number of columns of the completion above.
Let $\tilde h$ be the largest part of the strict completion
(Method~\ref{method:strict}).
It is computed from $v$ in $O(n^3)$ operations and, by
Lemma~\ref{lem:strict}, satisfies
\[
 \max(r_v,L)\leq\tilde h\leq\min(n-1,\ r_v+L-1).
\]
Put $\rho=\max(r_u,\tilde h)$, $\widetilde V=\fall{n}{\tilde h}$, and
define $R_q\leq2\rho$ and $U_q$ as in Section~\ref{sec:quantum}.

\begin{theorem}[Quantum, quantum double, and quantum triple bound]\label{thm:quantum}
For $v\neq\id$, Method~\ref{method:qproduct} computes every
coefficient $c_{u,v}^{w,\delta}(y;z)$ in \eqref{eq:qproduct}, as
arithmetic circuits over $\mathbb Z[q][y,z]$, in
\begin{equation}\label{eq:qmain}
 O\left(n^3+LM^3U_q^2\widetilde VQ_{\tilde h}(n)\right)
\end{equation}
elementary and coefficient-ring operations.
Every $w$ that occurs lies in $S_{n+L}$, and every $q^\delta$
that occurs involves only $q_1,\ldots,q_{n+L-1}$; hence the
coefficients at $z=y$ are the equivariant quantum
Littlewood--Richardson coefficients of $Fl_N$ for every $N\geq n+L$.
The bound is
\begin{equation}\label{eq:qcoarse}
 O\left(2^{2R_q}n^{2R_q+2\tilde h+4}\right)
 \subseteq O\left(2^{4\rho}n^{6\rho+4}\right)
 =2^{O(n\log n)},
\end{equation}
and it is polynomial in $n$ for fixed $r_u,r_v,L$.
The same bounds hold for quantum double ($z=y$) and quantum
($y=z=0$) products, again with the same state space.
\end{theorem}

Unlike Theorem~\ref{thm:main}, the polynomial case requires $L$ to be
fixed. For $\Sch_v=x_1^m$ one has $r_v=1$ but $\tilde h=m$;
see Section~\ref{sec:quantum}.

\subsection{Uniform bounds across six product types}
\label{sec:uniform}

Theorems~\ref{thm:main} and \ref{thm:quantum} cover ordinary, double,
triple, quantum, quantum double, and quantum triple products.
They are stated for the most general products, but the bounds are not
weakened by that generality.
The ordinary product is not analyzed separately and then lifted.
Instead, the triple computation and its ordinary specialization build
the same label graphs, make Pieri searches over the same chains, and
are bounded by the same state counts.
The two coefficient alphabets add only the cost of evaluating one
factorial elementary polynomial per matched pair
(Proposition~\ref{prop:uniform}).

This depends on two properties of the method.
\begin{enumerate}
\item The factorial elementary expansion of $\Sch_v(x;z)$
(Theorem~\ref{thm:elementary}) is indexed by the same filtered column
paths as the ordinary elementary expansion of $\Sch_v(x)$.
The alphabet $z$ enters only through the weight attached to each path.
\item The double Pieri rule \eqref{eq:pieri} is indexed by the same
canonical chains as Sottile's ordinary Pieri rule.
Each chain carries a single factorial elementary polynomial in a fixed
alphabet, rather than a sum over additional combinatorial data.
\end{enumerate}
The quantum Pieri rule \eqref{eq:qpieri} has the same feature, with
a $q$-monomial attached to each chain.
Its chains follow the quantum Pieri formula of Postnikov
\cite{postnikov1999quantum} and Ciocan-Fontanine
\cite{ciocan1999partial} rather than Sottile's rule.
Quantum covers add paths: they need not increase length, and the
strict completion can be taller than the ordinary one.
The extra paths are absorbed by the same counting, with $\tilde h$ and
$R_q$ in place of $h$ and $r$, so the asymptotic bound is essentially
unchanged (Theorem~\ref{thm:quantum}).

Neither property is a formal consequence of specialization.
Ordinary structure constants do determine the triple coefficients, but
the resulting reduction does not give a uniform bound.
By Hwang's coefficient lifting \cite[Theorem~1.2]{hwang2026lifting},
the coefficient of $y^\alpha$ in $c_{u,v}^w(y;z)$ is the single
coefficient $c_{\hat u,\hat v}^{\hat w}(0;z)$ for an explicit
block-lifted triple $(\hat u,\hat v,\hat w)$ in $S_{n+|\alpha|}$, where
$|\alpha|\leq\ell(u)+\ell(v)-\ell(w)$.
(Hwang expands in the Schubert basis of the second factor's alphabet,
so the two factors of \eqref{eq:product} are exchanged in his
notation; in ours the lifted alphabet $y$ is that of the output basis.)
The Cauchy formula then writes the coefficient of each monomial
$y^\alpha z^\beta$ as $\pm$ a sum of products of ordinary structure
constants \cite[Section~5.2]{hwang2026lifting}.
Thus an algorithm for $c_{u,v}^w$ does yield $c_{u,v}^w(y;z)$, but not
parsimoniously.
The reduction runs separately for each monomial and each Cauchy term.
Each run uses permutations whose rank exceeds $n$ by up to the degree of
the monomial.
Its cost is therefore multiplied by the expanded size of the triple
coefficients, and its input parameters grow.
Here, by contrast, the double and triple coefficients come at no extra
asymptotic cost because they are produced as arithmetic circuits along
the same paths as their ordinary specializations.
The circuit representation is essential to this statement, and we make
no comparable claim for expanded monomial output.

\subsection{Prior work and the encoding distinction}\label{sec:prior}

Pieri formulas for ordinary Schubert polynomials were established by
Sottile \cite{sottile}; we proved the double Pieri formula used here in
\cite[Theorem~7.1]{samuelmolev}.
The Cauchy identity and dominant factorization are standard
\cite{notes}.
We studied the triple coefficients in this generality in
\cite{samuelmolev}.
Hwang \cite{hwang2026lifting} recently showed that each of their
monomial coefficients is determined by ordinary structure constants
through an explicit block lift; Section~\ref{sec:uniform} compares that
reduction with the present method.
In earlier work \cite{samuelquantum}, we established the ordinary-double
expansion \cite[Theorem~4.5]{samuelquantum} and gave the corresponding
algorithms \cite[Methods~4.9--4.11]{samuelquantum}.
We also established the quantum and doubled Pieri formulas and the
quantum expansions \cite[Theorems~4.6--4.7]{samuelquantum}, which we use
in Section~\ref{sec:quantum}.
The contribution of the present paper is an explicit accounting of the
full-product algorithms in parameters obtainable from their inputs.

Pak and Slonim \cite[Theorem~1.4]{pakslonim} prove that a
\emph{single ordinary} coefficient $c_{u,v}^w$, for three permutations in
$S_n$ supplied by their range-restricted codes, can be computed in time
\begin{equation}\label{eq:pak}
 \exp\bigl(O(\psi^{5/2})\bigr)
       (\log n)^{O(\psi^2\log\psi)}.
\end{equation}
Here $\psi$ is their total-range parameter.
The bound is polynomial in $\log n$ for fixed $\psi$.
Their method expresses a coefficient as an alternating sum of
lattice-point counts and applies Barvinok's algorithm in bounded
dimension. The same polyhedral construction proves eventual
quasi-polynomiality of stretched Schubert coefficients
\cite[Theorem~1.1]{pakslonim}.

Our parameter $r$ measures how far the nonzero
Lehmer code extends from its first position.
We quote \eqref{eq:pak} with the parameter conventions of its source,
rather than deriving a formula relating $\psi$ to $r$.\footnote{The
version of \cite{pakslonim} consulted defines its range using exceedances
$w(i)>i$ and also states that the code vanishes off the exceedance set.
These two supports are not generally identical: $321$ has code $(2,1,0)$
and exceedance set $\{1\}$.
Our parameters and proofs use the Lehmer code directly and do not use
that identification.}

There are three important differences between \eqref{eq:pak} and
Theorem~\ref{thm:main}.
The former returns one integer.
The latter returns every coefficient of the triple product, each a
polynomial in two alphabets, within the same bound that it gives for
the ordinary integer coefficients.
The former is a bit-complexity statement for compressed input; the latter
initially counts arithmetic and elementary operations.
Finally, polynomial dependence on $n$ is much weaker than polynomial
dependence on $\log n$.
Thus the new bound does not subsume the bounded-range theorem.
At the unrestricted extreme our $2^{O(n\log n)}$ bound is smaller in
order than the bound obtained by substituting $\psi$ of order $n$
in \eqref{eq:pak}; this comparison of upper bounds is not a lower bound
on their algorithm or a claim about practical performance.

Encoding also matters for hardness statements.
Narayanan proves $\#P$-hardness for binary-encoded
Littlewood--Richardson coefficients \cite{narayanan2006complexity}.
These coefficients can be doubly exponential in the input length, so the
binary-encoded problem is not in $\#P$ in general.
For unary-encoded partitions, tableau counting places the problem in
$\#P$, while the cited work of Pak and Robichaux leaves its
$\#P$-hardness open \cite{pak2024vanishing}.
The associated Grassmannian Lehmer codes have a polynomial-size binary
encoding, so this gives $\#P$-hardness for the correspondingly compressed
Schubert problem.
It does not automatically give hardness for one-line permutations:
writing down a permutation associated with a large binary partition
part can require exponentially many entries.
Pak and Robichaux discuss this encoding distinction.
Their subsequent result places Schubert vanishing in probabilistic
polynomial time \cite{pakrobichaux2025bpp}, improving the earlier
polynomial-hierarchy result.
Neither vanishing nor a single-coefficient hardness statement is used
in our running-time proof.

\section{Dominant completion and a justified label filter}\label{sec:completion}

A permutation is \emph{dominant} if its code is a partition.
Its Rothe diagram is then a Young diagram.
If $\theta=\code(\tau^{-1})$ is its list of column heights, the
dominant factorization is
\begin{equation}\label{eq:dominant}
 \Sch_\tau(x;t)=\prod_{i=1}^{L}\prod_{j=1}^{\theta_i}(x_j-t_i).
\end{equation}

\begin{method}[Completion]\label{method:completion}
Start with $\mu=v^{-1}$.
While its code has a strict adjacent ascent
$\code_i(\mu)<\code_{i+1}(\mu)$, choose the largest such $i$
and replace $\mu$ by $\mu s_i$.
At termination put $\theta=\code(\mu)$ with trailing zeros deleted,
$\tau=\mu^{-1}$, and $a_*=v\mu$.
\end{method}

\begin{lemma}\label{lem:completion}
The completion remains in $S_n$ and terminates at a dominant permutation.
It satisfies
\[
 \ell(a_*)=\ell(\tau)-\ell(v),\qquad
 \Des(\tau)=\Des(v),\qquad
 \theta_1=r_v,\qquad L\leq n-1
\]
when $v\ne\id$.
It can be constructed in $O(n^3)$ elementary operations.
\end{lemma}
\begin{proof}
For an adjacent code ascent, right multiplication by $s_i$ changes
the two code entries $(a,b)$ to $(b+1,a)$ and increases length by one.
Thus there are at most $\binom n2$ steps and the final code is weakly
decreasing.
The successive right multiplications give a length-additive
factorization $\mu=v^{-1}a_*$, proving the length identity.

To check the descent assertion, consider the inverses before and after
one step. They differ by interchanging the consecutive values $i,i+1$.
Their positions are $\mu(i),\mu(i+1)$, with $\mu(i)<\mu(i+1)$.
These positions cannot be adjacent: if they were, the two code entries
of $\mu$ would be equal. Swapping consecutive values in nonadjacent
positions does not change any adjacent comparison.
The descent set of the inverse therefore stays unchanged throughout.
For a dominant permutation the number of nonzero rows of the diagram
is its last descent; transposing the diagram identifies this number
with $\theta_1$.
Finally, the diagram of any permutation in $S_n$ is contained in the
staircase, so $L\leq n-1$.
Maintaining the code, scanning for the next ascent, and updating the
permutation costs $O(n)$ per step.
\end{proof}

For $0\leq i\leq L$ define the dominant permutation $\tau_i$ by
\[
 \code(\tau_i^{-1})=(\theta_1,\ldots,\theta_i),\qquad
 F_i(x;t)=\Sch_{\tau_i}(x;t).
\]
Here $\tau_0=\id$, and
$F_i=F_{i-1}\prod_{j=1}^{\theta_i}(x_j-t_i)$.
All these permutations are in $S_n$.
Use the signed double-Schubert basis in the auxiliary variables
\[
 T_a(t;z)=(-1)^{\ell(a)}\Sch_a(t;z).
\]
We say $\sigma\leq_L\tau$ if $\tau=\rho\sigma$ with
$\ell(\tau)=\ell(\rho)+\ell(\sigma)$.

\begin{lemma}[Cauchy filter]\label{lem:filter}
Set
\[
 A_i=\{a\in S_n:\ell(a\tau_i)=\ell(\tau_i)-\ell(a)\}.
\]
Then
\begin{equation}\label{eq:filter}
 F_i(x;t)=\sum_{a\in A_i}\Sch_{a\tau_i}(x;z)T_a(t;z).
\end{equation}
Consequently $A_0=\{\id\}$, $a_*\in A_L$, and the coefficient of
$T_{a_*}$ in $F_L$ is $\Sch_v(x;z)$.
Moreover,
\begin{equation}\label{eq:labelcount}
 |A_i|\leq \fall{n}{h}=V.
\end{equation}
\end{lemma}
\begin{proof}
The double Cauchy identity and the interchange-of-alphabets identity give
\[
 \Sch_{\tau_i}(x;t)
 =\sum_{\substack{\rho\sigma=\tau_i\\ \ell(\rho)+\ell(\sigma)=\ell(\tau_i)}}
       \Sch_\sigma(x;z)\Sch_\rho(z;t),\qquad
 \Sch_\rho(z;t)=(-1)^{\ell(\rho)}\Sch_{\rho^{-1}}(t;z).
\]
Set $a=\rho^{-1}$. Then $\sigma=a\tau_i$ and these are exactly
the terms in \eqref{eq:filter}.
The claims about $a_*$ follow from Lemma~\ref{lem:completion}.

For the count, the map $a\mapsto a\tau_i$ is a bijection from $A_i$
onto the left weak interval $[\id,\tau_i]_L$.
If $\sigma\leq_L\tau_i$, every right descent of $\sigma$ is a
right descent of $\tau_i$: writing $\tau_i=\rho\sigma$,
a descent of $\sigma$ gives
$\ell(\tau_i s_j)\leq\ell(\rho)+\ell(\sigma)-1$.
The diagram of $\tau_i$ has at most $h$ rows, so these permutations
have no descent after position $h$.
A permutation in $S_n$ with an increasing suffix after position $h$
is determined by its ordered first $h$ values.
There are exactly $\fall{n}{h}$ such permutations.
\end{proof}

The filter is a polynomial identity, not an assertion that each discarded
path has zero weight. Coefficients of labels outside $A_i$ can vanish by
cancellation. They may be discarded \emph{at a column boundary} because
their entire coefficient is zero and subsequent multiplication is linear.
In particular, this does not justify pruning the interior of a column
merely because an intermediate path node fails the boundary test.

\section{Column expansions and the filtered path formula}\label{sec:columns}

For an ordered alphabet $Y=(Y_1,\ldots,Y_m)$ define
\begin{equation}\label{eq:elementary}
 E_{p,m}(Y;Z)=
 \sum_{1\leq j_1<\cdots<j_p\leq m}
       \prod_{q=1}^p (Y_{j_q}-Z_{j_q-q+1}),
 \qquad 0\leq p\leq m .
\end{equation}
Only the first $m-p+1$ entries of $Z$ are used when $p>0$.
Set $E_{0,m}=1$ and $E_{p,m}=0$ for $p>m$.
Repeated entries in $Z$ are allowed.

A \emph{pivot-$i$ path} of length $b$ is a sequence
$P=(a_0,\ldots,a_b)$ with
\[
 a_j=a_{j-1}t_{i,c_j},\qquad c_j\ne i,\qquad
 \ell(a_j)=\ell(a_{j-1})+1.
\]
No distinctness condition is imposed on the $c_j$.
Associate to it
\[
 \varepsilon(P)=(-1)^{\#\{j:c_j<i\}},\qquad
 Z(P)=(z_{a_0(i)},z_{a_1(i)},\ldots,z_{a_b(i)}).
\]

\begin{lemma}[Unrestricted column formula]\label{lem:column}
For $C_i=\prod_{j=1}^k(x_j-t_i)$,
\begin{equation}\label{eq:column}
 C_i T_a(t;z)=
 \sum_{b=0}^k\ \sum_{P:a_0=a}
   \varepsilon(P)E_{k-b,k}(x_1,\ldots,x_k;Z(P))T_{a_b}(t;z).
\end{equation}
\end{lemma}
\begin{proof}
The equivariant single-variable Monk formula, in the signed basis, is
\begin{equation}\label{eq:monk}
 (X-t_i)T_a=(X-z_{a(i)})T_a
   +\sum_{\substack{j>i\\ \ell(at_{ij})=\ell(a)+1}}T_{at_{ij}}
   -\sum_{\substack{j<i\\ \ell(at_{ji})=\ell(a)+1}}T_{at_{ji}} .
\end{equation}
It follows by subtracting the formulas for multiplication by
$t_1+\cdots+t_i$ and $t_1+\cdots+t_{i-1}$.
Apply \eqref{eq:monk} successively with $X=x_1,\ldots,x_k$.
Fix the subsequence of non-diagonal steps, namely a path $P$ of length $b$.
Its sign is $\varepsilon(P)$.
Choose the $k-b$ diagonal positions
$j_1<\cdots<j_{k-b}$ among the $k$ factors.
Before the diagonal position $j_q$ there have been $j_q-q$
non-diagonal steps, so its contribution is
$x_{j_q}-z_{a_{j_q-q}(i)}$.
Summing over the diagonal positions gives exactly
$E_{k-b,k}(x;Z(P))$ in \eqref{eq:elementary}.
\end{proof}

\begin{lemma}[Distinct-position cancellation]\label{lem:cancellation}
The sum in \eqref{eq:column} is unchanged if it is restricted to paths
whose nonpivot positions $c_1,\ldots,c_b$ are pairwise distinct.
For these paths $Z(P)$ consists of the values
\[
 D_i(a_0,a_b)=\{a_0(i)\}\cup
              \{a_0(j):a_0(j)\ne a_b(j)\},
\]
each once, and its order is immaterial.
\end{lemma}
\begin{proof}
Take a path with a repeated nonpivot position, and let the first
repetition occur at step $s$, with $c_s=c_t$ for $t<s$.
The positions $c_t,\ldots,c_{s-1}$ are distinct. Write them
$\gamma_1,\ldots,\gamma_m$, where $m=s-t$, and read indices cyclically.
For each $r$, consider the rotation
\[
 R_r=(\gamma_r,\gamma_{r+1},\ldots,\gamma_{r+m-1},\gamma_r).
\]
The rotation lemma \cite[Lemma~7.5]{samuelquantum} says that all valid
rotations have the same endpoint and, in cyclic order, alternate
between starting positions below and above the pivot $i$. Pair each
valid rotation starting above $i$ with the next valid one in cyclic
order, and each valid rotation starting below $i$ with the preceding
one. This is a fixed-point-free involution. Replacing the segment of
the original path by its paired rotation preserves the endpoints and
length, and the first repetition remains at step $s$.

The pivot visits the same multiset of values on the two paths: before
and after the segment the permutations agree, while during it the pivot
receives the values at $\gamma_1,\ldots,\gamma_m$ in a different order.
Thus the corresponding $E_{k-b,k}(x;Z)$ factors agree, by the symmetry
in the entries of $Z$ proved in \cite[Lemma~7.3]{samuelquantum}.
The two segments use each $\gamma_j$ once and repeat their respective
starting position once more. Those repeated positions lie on opposite
sides of $i$, so the signs $\varepsilon(P)$ are opposite. The terms
therefore cancel. Changing from the cited sign convention to
$\varepsilon(P)$ multiplies every length-$b$ term by $(-1)^b$; since
paired paths have the same length, this does not affect the cancellation.

On a distinct-position path the pivot successively carries
$a_0(i),a_0(c_1),\ldots,a_0(c_b)$.
These are exactly the values in $D_i(a_0,a_b)$.
The same symmetry gives the last assertion.
\end{proof}

\begin{theorem}[Filtered elementary expansion]\label{thm:elementary}
Start at label $\id$ and, for $i=1,\ldots,L$, take a pivot-$i$
path $P_i$ of length $b_i\leq\theta_i$ with pairwise distinct
nonpivot positions within that path.
Its starting label is the previous path's endpoint.
Require the endpoint at column $i$ to belong to $A_i$ and the final
endpoint to be $a_*$.
Then
\begin{equation}\label{eq:expansion}
 \Sch_v(x;z)=
 \sum_{(P_1,\ldots,P_L)}
 \prod_{i=1}^L
   \varepsilon(P_i)E_{\theta_i-b_i,\theta_i}(x;Z(P_i)).
\end{equation}
\end{theorem}
\begin{proof}
Apply Lemmas~\ref{lem:column} and \ref{lem:cancellation}
to the factors of $F_L$.
At each column boundary, Lemma~\ref{lem:filter} permits deleting all
labels outside $A_i$.
Extracting $T_{a_*}$ gives \eqref{eq:expansion}.
A path reaching $a_*\in S_n$ never leaves $S_n$, since all its
steps increase Bruhat order and the interval below a permutation in
$S_n$ is contained in $S_n$.
Before imposing the boundary filters, this is precisely
\cite[Theorem~4.5]{samuelquantum}.
\end{proof}

Setting $z=0$ in \eqref{eq:expansion} replaces $E_{p,m}$ by
$e_p(x_1,\ldots,x_m)$ and preserves the paths and their signs.
This is a statement about the \emph{path-indexed} expansion.
It is not a way of recovering the double polynomial from a collected
ordinary expansion after its path labels and alphabets have been erased.
No uniqueness or minimality of an elementary representation is needed.

\begin{method}[Backward construction of the label graphs]\label{method:labels}
Start with the singleton set $\{a_*\}$ at boundary $L$.
For $i=L,L-1,\ldots,1$ and each currently retained upper label,
enumerate all descending pivot-$i$ Bruhat paths of lengths
$0,\ldots,\theta_i$, never reusing a nonpivot position in a path.
Retain a path if its lower endpoint is in $A_{i-1}$.
Reverse it to obtain a forward edge, storing its two endpoints,
length, sign, and ordered alphabet indices.
Collect the distinct lower endpoints to seed the next boundary.
Do not prune internal path nodes using the endpoint filter.
\end{method}

This construction gives exactly the paths in Theorem~\ref{thm:elementary}
that can reach the required final label; it does not enumerate every
permutation in $A_i$.
Each stored edge represents an entire column path, not a single Bruhat cover.
The product of distinct pivot transpositions is a cycle, from which
their order can be recovered, so at most one such path has a specified
pair of endpoints.
This is the \texttt{kdown}/\texttt{vpaths} construction of
\cite[Methods~4.9--4.10]{samuelquantum}.
In its implementation, a nonpivot position is eligible exactly when it
still holds its value from the start of the path: distinct pivot swaps
change every position already used and leave all others untouched.

\section{Pieri multiplication and coefficient construction}\label{sec:pieri}

For a permutation $w$ and an integer $k$, a canonical elementary
Pieri chain consists of Bruhat covers
$w_j=w_{j-1}t_{a_j b_j}$ with $a_j\leq k<b_j$,
whose moved-in values $w_{j-1}(b_j)$ are strictly decreasing.
The elementary Pieri rule identifies chains of length $a$ with the
terms of $e_a(x_1,\ldots,x_k)\Sch_w(x)$, each with multiplicity one
\cite{sottile}.
In particular, distinct canonical chains from a fixed $w$ have distinct
endpoints, including across different lengths.
No left position can be used twice: a second ascent there would require
a moved-in value greater than the first, contrary to the ordering.

Write $w\longrightarrow_k w'$ for such an endpoint and let
$a=\ell(w')-\ell(w)$.
The fixed alphabet is
\[
 Y(w,w';k)=(y_{w(j)}:1\leq j\leq k,\ w(j)=w'(j)),
\]
listed in increasing order of $j$; it has $k-a$ entries.
The double Pieri formula \cite[Theorem~7.1]{samuelmolev} is
\begin{equation}\label{eq:pieri}
 \Sch_w(x;y)E_{p,k}(x;Z)
 =\sum_{\substack{w\longrightarrow_k w'\\ a\leq p}}
     E_{p-a,k-a}(Y(w,w';k);Z)\Sch_{w'}(x;y).
\end{equation}
Since this is a polynomial identity in the coefficient alphabet,
it remains valid when entries of $Z$ coincide, as they may in
Lemma~\ref{lem:column}.

\begin{method}[Pieri enumeration]\label{method:pieri}
Maintain a frontier $(w,\beta)$, initially $(w,+\infty)$.
To extend it, test $a\leq k<b$ for
\[
 w(a)<w(b)<\beta,\qquad
 \nexists c\ (a<c<b,\ w(a)<w(c)<w(b)).
\]
For each successful test append $(wt_{ab},w(b))$ to the next frontier.
Return the initial permutation and every endpoint through the requested
maximum degree $p\leq k$, recording its degree.
\end{method}

\begin{proposition}\label{prop:pieriwork}
If all generated permutations lie in $S_M$ and $B$ endpoints are
returned, Method~\ref{method:pieri} costs $O(M^3B)$ elementary operations.
\end{proposition}
\begin{proof}
A frontier node has at most $k(M-k)$ candidate transpositions.
Testing the Bruhat-cover condition and copying a permutation each cost
$O(M)$.
Multiplicity-freeness of the canonical Pieri chains bounds the total
number of processed nodes by the number of returned endpoints.
\end{proof}

\begin{lemma}[Polynomial-cost elementary weights]\label{lem:weight}
A factor $E_{p,m}(Y;Z)$ can be built as an arithmetic circuit using
$O((m+1)^2)$ ring operations.
\end{lemma}
\begin{proof}
Separating subsets according to whether they contain the last index gives
\begin{equation}\label{eq:weightrec}
 E_{p,m}(Y;Z)=E_{p,m-1}(Y_1,\ldots,Y_{m-1};Z)
 +(Y_m-Z_{m-p+1})
     E_{p-1,m-1}(Y_1,\ldots,Y_{m-1};Z).
\end{equation}
Use $E_{0,t}=1$ and $E_{q,t}=0$ for $q>t$, and memoize the
subproblems $(q,t)$ reached from $(p,m)$.
There are at most $(m+1)^2$ of them, each costing a bounded number
of operations.
Every nontrivial reached subproblem has $t-q\leq m-p$,
so no coefficient variable beyond $Z_{m-p+1}$ is requested.
The circuit shares subexpressions rather than copying expression trees.
\end{proof}

An unmemoized divide-and-conquer evaluation is another valid choice:
split $Y$ after $s=\lfloor m/2\rfloor$ and use
\[
 E_{p,m}(Y;Z)=\sum_a E_{a,s}(Y';Z)
     E_{p-a,m-s}(Y'';Z^{[s-a]}),\qquad
 Z^{[c]}_j=Z_{j+c},
\]
where the sum is over feasible degrees.
Its recursion gives an upper bound $m^{O(\log m)}$.
That is an evaluation-strategy cost, not an intrinsic lower bound
for factorial elementary polynomials.
Theorem~\ref{thm:main} uses \eqref{eq:weightrec}.

\begin{method}[Classical \alg]\label{method:product}
Construct the completion and the label graphs.
The coefficient ring is $\mathbb Z[y,z]$.
Initialize a table $H$ with $H(u,\id)=1$.
At column $i$ put $k=\theta_i$ and create a fresh table $H'$.
\begin{enumerate}
\item Group the keys of $H$ by their first component $w$.
For each distinct $w$ enumerate its Pieri endpoints once, through degree
$\min(k,d-\ell(w))$.
Reuse this list for all labels paired with $w$.
\item For each table entry $H(w,a)$, each label edge $P:a\to a'$
of length $b$, and each Pieri endpoint $w'$ of degree $a_u$,
use the pair only if $a_u+b\leq k$.
Accumulate
\[
 H'(w',a')\mathrel{+}=
  H(w,a)\varepsilon(P)
  E_{k-a_u-b,k-a_u}\bigl(Y(w,w';k);Z(P)\bigr).
\]
\item Replace $H$ by $H'$ and continue.
After column $L$, return $H(w,a_*)$ as the coefficient of $\Sch_w(x;y)$.
\end{enumerate}
In the ordinary specialization keep only $a_u+b=k$;
the elementary weight is then $1$.
\end{method}

\begin{theorem}[Correctness]\label{thm:correct}
Method~\ref{method:product} computes \eqref{eq:product}.
\end{theorem}
\begin{proof}
Multiply the expansion in Theorem~\ref{thm:elementary} by
$\Sch_u(x;y)$ and apply \eqref{eq:pieri} at each column.
The displayed update is exactly the product of the column sign,
the running coefficient, and the Pieri coefficient.
Dynamic programming merely combines contributions with equal endpoint
pairs.
The final label condition is precisely that in
Theorem~\ref{thm:elementary}.

It remains to justify the degree cap.
Along a complete label path,
$\sum_i b_i=\ell(a_*)=\ell(\tau)-\ell(v)$ and
$\sum_i\theta_i=\ell(\tau)$.
Every admissible paired path has
$\sum_i a_{u,i}\leq\sum_i(\theta_i-b_i)=\ell(v)$.
All $a_{u,i}$ are nonnegative, so no contributing prefix can reach
a permutation of length greater than $d$.
Thus the cap discards no contribution.
Specialization to $y=z=0$ annihilates every elementary weight of
positive degree and leaves precisely the stated ordinary update.
\end{proof}

\section{Input bounds for every enumerated state}\label{sec:states}

\begin{lemma}[Pieri ambient size]\label{lem:ambient}
If $w\in S_m$ and $k<m$, a canonical elementary Pieri chain from
$w$ stays in $S_{m+1}$.
\end{lemma}
\begin{proof}
Before any position beyond $m$ is used, a Bruhat cover cannot have
right endpoint $b\geq m+2$: position $m+1$ has value $m+1$,
strictly between $w(a)\leq m$ and $w(b)=b$.
The only possible new value is therefore $m+1$.
If it is moved in, it must be the first moved-in value, since all
previously available values are at most $m$ and the sequence decreases.
All subsequent moved-in values are smaller; no larger new position
can be used.
\end{proof}

\begin{lemma}[Code-length and degree bounds]\label{lem:states}
Every $u$-side permutation enumerated by Method~\ref{method:product},
including nodes of Pieri searches that do not survive into the next
table, lies in $S_M$, has no descent after $r$, and has length at most $d$.
There are at most $U$ such permutations.
\end{lemma}
\begin{proof}
The ambient claim follows by applying Lemma~\ref{lem:ambient} once
per column, starting in $S_n$.
The degree cap in Method~\ref{method:product} bounds every partial
Pieri chain by $d$.

For the descent claim, let $w$ have no descent after $r$.
Since $k\leq h\leq r$, the polynomial
\[
 e_a(x_1,\ldots,x_k)\Sch_w(x)
\]
involves only $x_1,\ldots,x_r$.
For $j>r$ its divided difference $\partial_j$ is zero.
For a Schubert basis element, $\partial_j\Sch_b=\Sch_{bs_j}$
if $j\in\Des(b)$; otherwise its divided difference is zero.
Distinct $b$ with descent $j$ give distinct $bs_j$.
Linear independence of the Schubert basis therefore implies that no
permutation in this product has a descent after $r$.
Every canonical chain of length $a$ is a term of this ordinary Pieri
product with coefficient one, so this argument applies to
every enumerated endpoint, before any double-coefficient cancellation.
Induction starts with $u$.

Permutations of $S_M$ with an increasing suffix after $r$ number
$\fall{M}{r}$.
Alternatively, their Lehmer codes have length at most $r$ and sum
at most $d$. There are $\binom{d+r}{r}$ nonnegative $r$-tuples with
that sum bound, by stars and bars.
Taking the smaller estimate proves the assertion.
\end{proof}

\begin{lemma}[Cost of the label construction]\label{lem:labelwork}
Method~\ref{method:labels} stores at most $L V Q_h(n)$ column paths
and costs $O(n^3+L n^3 V Q_h(n))$ elementary operations.
This estimate charges paths rejected by the endpoint filter.
\end{lemma}
\begin{proof}
There are at most $V$ retained upper labels at a boundary, by
Lemma~\ref{lem:filter}.
From each upper label, a descending pivot path has depth at most
$\theta_i\leq h$ and uses each nonpivot position at most once.
There are at most $\fall{n-1}{b}$ possible position sequences of
length $b$, even before testing the Bruhat conditions.
The entire search tree, including rejected nodes, consequently has
at most $Q_h(n)$ nodes.
Every node is in $S_n$.
Scanning candidates, testing Bruhat covers, computing the endpoint
length test, and recording a path take polynomial work bounded by
$O(n^3)$ per node.
Summing over boundaries and upper labels proves the claim.
Constructing the prefixes $\tau_i$ and the completion fits the
additional $O(n^3)$ term.
\end{proof}

\section{Complexity and its interpretation}\label{sec:complexity}

\subsection{Arithmetic model}

A coefficient is stored as a directed acyclic arithmetic circuit over
$\mathbb Z[y,z]$; addition and multiplication create gates pointing to
existing gates.
The elementary routine shares its dynamic-programming subexpressions.
We count these ring operations and elementary operations on permutation
entries separately from bit costs.
Permutation entries have $O(\log M)$ bits.
Deterministic dictionaries may be implemented by balanced trees:
comparing two keys costs $O(M)$ elementary operations, and a table
has at most $UV$ keys.
Since $\log(UV)=O(M\log M)$, this bookkeeping is absorbed by the
$M^3$ factor below.
No unit-cost polynomial identity test is assumed.

The returned table is an exact expansion but need not be a
minimal-support list.
An identically zero coefficient can be left as a circuit computing zero.
Removing all such entries, expanding coefficients into monomials, and
canonical symbolic simplification are different output requirements
and are not included in the theorem.

\begin{proof}[Proof of Theorem~\ref{thm:main}]
The label preprocessing costs
$O(n^3+L n^3 VQ_h(n))$ by Lemma~\ref{lem:labelwork}.
At a product boundary there are at most $U$ distinct $u$-side
permutations and at most $V$ labels.
At one column the Pieri enumeration is called at most $U$ times,
each returning at most $U$ endpoints by Lemma~\ref{lem:states}.
Proposition~\ref{prop:pieriwork} bounds this work by $O(M^3U^2)$.

Count matched pairs by iterating over all column edges, rather than
assuming a bound on the forward outdegree of each label.
There are at most $VQ_h(n)$ such edges in a column.
For each there are at most $U$ possible source permutations,
and each source has at most $U$ Pieri endpoints.
Thus there are at most $U^2VQ_h(n)$ pairs.
Forming the fixed alphabet takes $O(M)$ work.
Lemma~\ref{lem:weight} builds its elementary coefficient in
$O((h+1)^2)$ ring operations; the remaining update takes a bounded
number of ring operations.
Permutation and dictionary work, as well as this arithmetic, is
bounded by $O(M^3)$ per pair.
Summing over $L$ columns and absorbing the preprocessing and
Pieri terms gives \eqref{eq:main}.

For the coarse estimate use
\[
 L\leq n-1,\quad M\leq2n-1,\quad U\leq(2n)^r,\quad
 V\leq n^h,\quad Q_h(n)\leq2n^h .
 \]
The last bound follows by comparing $\fall{n-1}{b}$ with $(n-1)^b$
and summing the geometric series for $n\geq3$, and
for $n=2$ because $h\leq1$.
Substitution gives \eqref{eq:coarse}.
Since $h\leq r\leq n-1$, its logarithm is at most
$4n\log_2 n+O(n)$.

If $v$ is dominant, $\mu=v^{-1}$ and $a_*=\id$.
No descending path of positive length starts at the identity.
Every label graph therefore consists of the single degree-zero edge
$\id\to\id$, and the entire product calculation has one label.
Replacing $VQ_h(n)$ by $1$ gives the dominant-input estimate.
By Proposition~\ref{prop:uniform}, the double and ordinary products
satisfy the same bound.
\end{proof}

\subsection{Uniformity across coefficient rings}

The ordinary specialization of Method~\ref{method:product} keeps only
pairs with $a_u+b=k$. In the double specialization, set $z=y$.

\begin{proposition}[Uniformity]\label{prop:uniform}
Fix $u,v$ with $v\ne\id$.
Compare the triple computation of Method~\ref{method:product} with its
double and ordinary specializations.
\begin{enumerate}
\item All three use the same completion and the same label graphs.
\item Every table key and every Pieri call made by the ordinary
computation is made by the triple computation, with the same arguments.
Every matched pair used by the ordinary computation is used by the
triple computation.
The triple computation's additional pairs are those with $a_u+b<k$.
\item Per matched pair, the triple and double computations perform
$O((h+1)^2)$ ring operations and $O(M)$ elementary operations beyond
the ordinary computation.
\end{enumerate}
In particular, all three computations satisfy \eqref{eq:main}.
The same statements compare the quantum triple computation of
Method~\ref{method:qproduct} with its quantum double and quantum
specializations, with
$t+b=k$ in place of $a_u+b=k$ and $\tilde h$ in place of $h$, except
that a $\mathtt{qpieri}$ call of the quantum triple computation may have a
larger depth $p_w$. Its returned list then contains that of the
quantum call.
\end{proposition}
\begin{proof}
The completion and the label graphs depend only on $v$.
The recorded alphabet indices $Z(P)$ are used only to form weights,
never to accept or reject a path.
Step~(1) of Method~\ref{method:product} depends only on $w$, $k$, and $d$.
Step~(2) accepts a pair according to $a_u$ and $b$ alone.
By induction over columns, the ordinary keys at each boundary therefore
form a subset of the triple keys, and the ordinary pairs form a subset
of the triple pairs.
The additional work for a triple pair is the formation of $Y(w,w';k)$
and the circuit of Lemma~\ref{lem:weight}.
In the ordinary computation the weight is $1$.
The proofs of Lemmas~\ref{lem:filter}, \ref{lem:states}, and
\ref{lem:labelwork} never use the coefficient ring, so they apply to all
three computations.
For the quantum triple method, the same argument applies with
Lemmas~\ref{lem:strict} and \ref{lem:qstates}.
A $\mathtt{qpieri}$ search of depth $p$ is a prefix of one of larger
depth, and Lemma~\ref{lem:qstates} bounds every chain that the
quantum triple computation enumerates.
\end{proof}

The proposition is an upper-bound statement.
The ordinary computation can be cheaper in practice, since it discards
pairs with $a_u+b<k$.
The point is that the worst-case analysis does not distinguish the
three problems, and the bound for the triple product is not obtained by
paying for its extra generality.

\begin{corollary}[Ordinary bit costs]\label{cor:bits}
For unit initial coefficient, the ordinary algorithm has
bit complexity bounded by \eqref{eq:main} times a polynomial in $n$.
It therefore still has complexity $2^{O(n\log n)}$, and polynomial
complexity in $n$ for fixed $r_u,r_v$.
\end{corollary}
\begin{proof}
Each surviving matched pair contributes the running integer with
sign $+1$ or $-1$.
The absolute value of any intermediate sum is bounded by the number
of complete or partial paired paths.
A state has at most $UVQ_h(n)$ successors counted with multiplicity,
so that number is at most $(2UVQ_h(n))^L$.
Its bit length is
\[
 O\bigl(1+L\log(2UVQ_h(n))\bigr)=O(n^2\log n).
 \]
Integer additions and permutation operations consequently incur only
a polynomial bit overhead.
This argument bounds intermediate signed sums as well as final
nonnegative coefficients; it does not rely on positivity preventing
intermediate cancellation.
\end{proof}

\subsection{Dependence on code length}

The number of code entries is a useful parameter because the
$u$-side count can be replaced by $\binom{d+r}{r}$ when degree
is more restrictive than the ambient permutation size.
No weak-order interval enumeration is needed for this refinement.
The $v$-side bound depends on $h=r_v$, reflecting the fact that the
second factor supplies the elementary expansion.
For ordinary multiplication the roles may be interchanged and the
smaller of the two bounds used.
For \eqref{eq:product}, such an interchange also changes which
coefficient alphabet indexes the output basis and is not automatic.

An input code with $r_w$ entries has encoding length
$O(r_w\log n)$, whereas the algorithm explicitly works with
permutations and has polynomial dependence on $n$ for fixed $r$.
This is a pseudo-polynomial feature relative to binary code entries.
For example, $\Sch_{\operatorname{uncode}(0,A)}(x_1,x_2)=h_A(x_1,x_2)$,
and
\[
 h_A(x_1,x_2)^2
   =\sum_{j=0}^{A}s_{(2A-j,j)}(x_1,x_2).
\]
The two input codes have constant length and $O(\log(A+1))$ bits,
but an explicit Schubert-basis list has $A+1$ terms.

Under scaling of the entries of fixed Lehmer codes by a factor $N$,
the code lengths remain fixed and $n=O(N)$.
Theorem~\ref{thm:main} hence gives a polynomial-in-$N$ bound for these
full products.
This does not yield the polylogarithmic dependence available from
\eqref{eq:pak} in its fixed-range regime, nor does it prove eventual
quasi-polynomiality.

\begin{corollary}[Grassmannian inputs]\label{cor:grass}
Let $u,v\in S_n$, $v\ne\id$, have no descents other than at $k$.
Then every column of the completion has height $k$, and every $u$-side
permutation enumerated by Method~\ref{method:product}, including
Pieri search nodes, is $k$-Grassmannian.
Every label $a\in A_i$ has $a\tau_i$ $k$-Grassmannian.
Consequently $U$ and $V$ in \eqref{eq:main} may be replaced by
\[
 U_k=\min\left\{\binom{M}{k},\binom{d+k}{k}\right\},\qquad
 V_k=\binom nk ,
\]
giving
\[
 O\left(n^3+LM^3U_k^2V_kQ_k(n)\right)
 =O\left(n^{4}\binom{2n}{k}^{2}\binom nk n^{k}\right).
\]
Uniformly in $k$ this is $O(4^kn^{4k+4}/(k!)^3)$, and in general it is
$2^{O(n)}n^{k+4}$.
\end{corollary}
\begin{proof}
By Lemma~\ref{lem:completion}, $\Des(\tau)=\Des(v)=\{k\}$.
A dominant permutation whose only descent is at $k$ has a rectangular
diagram with $k$ rows, so its transpose $\theta$ has all parts equal to
$k$.
For the labels, $a\tau_i\leq_L\tau_i$, so as in the proof of
Lemma~\ref{lem:filter} its descents lie in $\Des(\tau_i)=\{k\}$;
there are at most $\binom nk$ such permutations in $S_n$.

Let $w$ be $k$-Grassmannian. Then $\Sch_w(x)$ is symmetric in
$x_1,\ldots,x_k$ and involves no other $x$-variables, so the same
holds for $e_a(x_1,\ldots,x_k)\Sch_w(x)$.
Hence $\partial_j$ annihilates it for every $j\ne k$.
The argument of Lemma~\ref{lem:states} shows that every permutation in
this ordinary Pieri product, and so every canonical chain endpoint,
has no descent other than at $k$.
Induction from $u$ covers all enumerated states, including those of
the double computation, whose Pieri chains are the same.
A $k$-Grassmannian permutation in $S_M$ is determined by its first $k$
values, as a set, giving $\binom Mk$.
Alternatively, its code is a partition with at most $k$ parts,
read in reverse, and of size at most $d$; there are at most
$\binom{d+k}{k}$ of these.
The displayed estimate follows from the proof of Theorem~\ref{thm:main},
with $M\leq 2n$ and $Q_k(n)\leq 2n^k$.
\end{proof}

Thus Grassmannian multiplication works entirely with partitions in a
$k\times(M-k)$ box on the $u$-side, with
$L=\code_k(v)\leq n-k$ columns.
The remaining non-exponential factor $n^{k}$ comes from the label
search $Q_k(n)$, which counts rejected distinct-position paths.
Sharpening that count would make the Grassmannian bound $2^{O(n)}$.

\section{Quantum, quantum double, and quantum triple multiplication}\label{sec:quantum}

We now prove Theorem~\ref{thm:quantum}, directly for the quantum triple
product.
The architecture is unchanged: a filtered elementary expansion of the
second factor is paired with a Pieri rule for the first.
Two ingredients change.
The factors must be quantized one column at a time, which forces
pairwise distinct column heights; and the quantum Pieri rule has
length-decreasing steps, which requires new state estimates.

\subsection{Quantization and the strict completion}

Let $X_j=(x_1,\ldots,x_j)$.
A \emph{standard elementary monomial} is a product
$\prod_j e_{c_j}(X_j)$ in which each height $j$ occurs at most once and
$0\leq c_j\leq j$.
These monomials form a $\mathbb Z$-basis of $\mathbb Z[x]$.
The quantization map of Fomin, Gelfand, and Postnikov
\cite{fomin1997quantum} is the $\mathbb Z[q]$-linear map, extended
linearly over all coefficient alphabets, sending
$\prod_j e_{c_j}(X_j)$ to $\prod_j e^q_{c_j}(X_j)$, where
$e^q_c(X_j)$ are the quantum elementary polynomials.
It sends $\Sch_w(x;y)$ to the quantum double Schubert polynomial
$\widetilde\Sch_w(x;y)$ \cite{fomin1997quantum,kirillov2000quantum},
\cite[Lemma~2.1]{samuelquantum}.
These polynomials, indexed by $S_\infty$, form a
$\mathbb Z[q][y]$-basis of $\mathbb Z[q][x,y]$, so we may write
\begin{equation}\label{eq:qproduct}
 \widetilde\Sch_u(x;y)\widetilde\Sch_v(x;z)
  =\sum_{w,\delta}q^\delta c_{u,v}^{w,\delta}(y;z)\,
     \widetilde\Sch_w(x;y).
\end{equation}
Equation~\eqref{eq:qproduct} is the quantum triple product.
At $z=y$ it becomes the quantum double product, whose coefficient
$c_{u,v}^{w,\delta}(y;y)$ is the
equivariant quantum Littlewood--Richardson coefficient of $Fl_N$ for every
$N$ such that every occurring $w$ lies in $S_N$ and every occurring degree
has $\delta_i=0$ for $i\geq N$
\cite{kim1996equivariant,lam2014quantum,mihalcea2006equivariant},
\cite[Section~2]{samuelquantum}.
At $y=z=0$ it becomes the quantum Schubert product, whose coefficients
are the non-equivariant quantum structure constants.
Grading $\deg q_i=2$ makes \eqref{eq:qproduct} homogeneous.
Put
\[
 \widetilde E_{p,k}(x;Z)=\mathcal Q\bigl(E_{p,k}(x;Z)\bigr),
\]
where $\mathcal Q$ is the quantization map.
Expanding the definition, $E_{p,k}(x;Z)$ is a $\mathbb Z[Z]$-linear
combination of $e_c(X_k)$ with $c\leq k$, so $\widetilde E_{p,k}$ is
the corresponding combination of $e^q_c(X_k)$
\cite[Lemma~7.3]{samuelquantum}.

The classical completion of Section~\ref{sec:completion} usually has
repeated column heights, and $\mathcal Q$ is not multiplicative on
$e_a(X_k)e_b(X_k)$.
We therefore use the strict completion of
\cite[Definition~4.1]{samuelquantum}.
Call $\mu$ \emph{strict dominant} if it is dominant and its nonzero
code entries are strictly decreasing.

\begin{method}[Strict completion]\label{method:strict}
Start with $\mu=v^{-1}$.
While $\mu$ is not strict dominant, choose the largest $i$ with
$\code_{i+1}(\mu)\neq0$ and $\code_i(\mu)\leq\code_{i+1}(\mu)$,
and replace $\mu$ by $\mu s_i$.
At termination put $\tilde\theta=\code(\mu)$ with trailing zeros
deleted, $\tilde\tau=\mu^{-1}$, $\tilde a_*=v\mu$, and
$\tilde h=\tilde\theta_1$.
\end{method}

\begin{lemma}[Strict height]\label{lem:strict}
Let $v\ne\id$ and let $L$ be the number of nonzero entries of
$\code(v^{-1})$.
The strict completion stays in $S_n$, uses $O(n^3)$ elementary
operations, and satisfies
$\ell(\tilde a_*)=\ell(\tilde\tau)-\ell(v)$.
Its shape has exactly $L$ parts,
$\tilde\theta_1>\cdots>\tilde\theta_L>0$,
and the classical completion also has exactly $L$ columns.
Finally,
\begin{equation}\label{eq:strictheight}
 \max(r_v,L)\leq\tilde h\leq
 \min\bigl(n-1,\ r_v+L-1\bigr).
\end{equation}
\end{lemma}
\begin{proof}
A step at $i$ has $\code_i(\mu)\leq\code_{i+1}(\mu)$, so
$\mu(i)<\mu(i+1)$; it replaces the entries $(a,b)$ by $(b+1,a)$ and
increases length by one.
Since $b\neq0$, we have $i+1\leq m-1$ whenever $\mu\in S_m$, so
$\mu s_i\in S_m$.
This proves the ambient claim for $m=n$, termination, the length identity,
and the $O(n^3)$ cost as in Lemma~\ref{lem:completion}.
The replacement $(a,b)\mapsto(b+1,a)$ with $b\neq0$ preserves the
number of nonzero entries. The same holds for the classical steps,
which have $a<b$.
Hence both shapes have $L$ parts, and the termination condition gives
strictly decreasing parts.

A strictly decreasing sequence of $L$ positive integers has first entry at
least $L$.
For the other lower bound, $\tilde\tau=\tilde a_*^{-1}v$ is
length-additive,
so, as in the proof of Lemma~\ref{lem:filter}, every right descent of
$v$ is one of $\tilde\tau$.
The permutation $\tilde\tau$ is the inverse of a dominant permutation,
hence dominant, and its last descent is the number of rows of its
diagram, namely $\tilde\theta_1$.
Thus $r_v\leq\tilde h$.

For the upper bound let $\mu$ be the classical completion, with shape
$\theta$, and put $m=\max_j(\theta_j+j)$.
A dominant permutation with shape $\theta$ lies in $S_m$:
its code satisfies $\code_j\leq m-j$.
Since $\mu=v^{-1}a_*$ with lengths adding, $v^{-1}$ is a prefix of a
reduced word of $\mu$ and so lies in $S_m$.
Applying the first paragraph with this $m$ shows that the strict
completion lies in $S_m$, whence $\tilde h=\code_1\leq m-1$.
Finally, $\theta_j\leq\theta_1=r_v$ and $j\leq L$ give
$m\leq r_v+L$. The bound $\tilde h\leq n-1$ is the case $m=n$.
\end{proof}

Both inequalities in \eqref{eq:strictheight} can be sharp.
For $v=(m+1)12\cdots m$, which has $\Sch_v=x_1^m$, one has $r_v=1$,
$L=m$, $\theta=(1^m)$, and $\tilde\theta=(m,m-1,\ldots,1)$.
Computations through $S_8$ also gave
$\tilde h=\max_j(\theta_j+j)-1$ for every nonidentity $v$.

\subsection{Quantum filtered expansion}

Define $\tilde\tau_i$, $\tilde F_i$, $\tilde A_i$, and the label graphs
exactly as in Sections~\ref{sec:completion} and \ref{sec:columns},
with $\tilde\theta$ in place of $\theta$, and put
\[
 \widetilde V=\fall{n}{\tilde h}.
\]
The proofs of Lemmas~\ref{lem:filter} and \ref{lem:labelwork} use only
that the completion is dominant, lies in $S_n$, and satisfies
the length identity, and that its first column has height
equal to the number of rows of the diagram of each $\tilde\tau_i$.
They therefore give $|\tilde A_i|\leq\widetilde V$, and a label
construction with at most $L\widetilde VQ_{\tilde h}(n)$ edges
costing $O(n^3+Ln^3\widetilde VQ_{\tilde h}(n))$.
Likewise, the proof of Theorem~\ref{thm:elementary} gives
\begin{equation}\label{eq:strictexpansion}
 \Sch_v(x;z)=
 \sum_{(P_1,\ldots,P_L)}
 \prod_{i=1}^L
   \varepsilon(P_i)E_{\tilde\theta_i-b_i,\tilde\theta_i}(x;Z(P_i)),
\end{equation}
summed over the filtered distinct-position paths ending at
$\tilde a_*$.

\begin{theorem}[Filtered quantum expansion]\label{thm:qexpansion}
With the same index set as \eqref{eq:strictexpansion},
\begin{equation}\label{eq:qexpansion}
 \widetilde\Sch_v(x;z)=
 \sum_{(P_1,\ldots,P_L)}
 \prod_{i=1}^L
   \varepsilon(P_i)\widetilde E_{\tilde\theta_i-b_i,\tilde\theta_i}
     (x;Z(P_i)).
\end{equation}
\end{theorem}
\begin{proof}
Apply $\mathcal Q$ to \eqref{eq:strictexpansion}.
On the left this gives $\widetilde\Sch_v(x;z)$.
Each factor in a summand is a $\mathbb Z[z]$-combination of
$e_c(X_{\tilde\theta_i})$, and the heights $\tilde\theta_i$ are
pairwise distinct.
Expanding a summand therefore gives a $\mathbb Z[z]$-combination of
standard elementary monomials.
Hence $\mathcal Q$ maps the summand to the product of the quantized
factors.
\end{proof}

Without the filters, \eqref{eq:qexpansion} is
\cite[Theorem~4.6]{samuelquantum}.
The proof shows more: for \emph{any} set of complete distinct-height
column paths, the sum of quantized weights is the quantization of
the sum of classical weights.
Thus every boundary filter justified classically is automatically
justified in quantum cohomology, without a new cancellation argument.

\subsection{The quantum Pieri step}

A \emph{quantum cover} of $w$ is a transposition $t_{ab}$, $a<b$, with
$\ell(wt_{ab})=\ell(w)+1$, or with $w(a)>w(b)$ and
$\ell(wt_{ab})=\ell(w)+1-2(b-a)$.
The latter is equivalent to every $w(c)$, $a<c<b$, lying strictly
between $w(b)$ and $w(a)$, and it carries the weight
$q_a\cdots q_{b-1}$.
The routine $\mathtt{qpieri}(w,p,k)$ of
\cite[Method~4.12]{samuelquantum} is Method~\ref{method:pieri} with
quantum covers in place of Bruhat covers, with $a\leq k<b$, with each
left position used at most once, and with the right positions $b$
weakly decreasing along a chain.
It returns triples $(w',t,q^{\delta'})$, one for every node of its
search tree.
By \cite[Lemma~4.13]{samuelquantum} it returns each endpoint once and,
for $p'\leq p\leq k$,
\begin{equation}\label{eq:qpieri}
 \widetilde\Sch_w(x;y)\widetilde E_{p',k}(x;Z)
 =\sum_{\substack{(w',t,q^{\delta'})\\ t\leq p'}}
   q^{\delta'}E_{p'-t,k-t}\bigl(Y(w,w';k);Z\bigr)
   \widetilde\Sch_{w'}(x;y).
\end{equation}
Here $t$ is the number of transpositions in the chain, and
$Y(w,w';k)$ is as in Section~\ref{sec:pieri}.
Since left positions are used once, it has $k-t$ entries.
Multiplicity-freeness gives the analogue of
Proposition~\ref{prop:pieriwork}: if all nodes lie in $S_M$ and $B$
triples are returned, the call costs $O(M^3B)$ elementary operations.
The $q$-weight can be maintained by range-increment updates to its
exponent vector and materialized once per returned chain, for
$O(MB)$ additional operations.

\begin{method}[Quantum triple \alg]\label{method:qproduct}
Construct the strict completion and its label graphs.
The coefficient ring is $\mathbb Z[q][y,z]$.
Initialize $H(u,\id)=1$.
At column $i$ put $k=\tilde\theta_i$ and create a fresh table $H'$.
\begin{enumerate}
\item Group the keys of $H$ by their first component $w$.
For each such $w$ let $p_w$ be the maximum of $k-b$ over all label
edges, of length $b$, leaving labels paired with $w$.
Call $\mathtt{qpieri}(w,p_w,k)$ once.
\item For each table entry $H(w,a)$, each edge $P:a\to a'$ of length $b$,
and each returned triple $(w',t,q^{\delta'})$ with $t+b\leq k$, accumulate
\[
 H'(w',a')\mathrel{+}=
  H(w,a)\,\varepsilon(P)\,q^{\delta'}\,
  E_{k-t-b,k-t}\bigl(Y(w,w';k);Z(P)\bigr).
\]
\item Replace $H$ by $H'$.
After column $L$, return $H(w,\tilde a_*)$; its coefficient of
$q^\delta$ is $c_{u,v}^{w,\delta}(y;z)$.
\end{enumerate}
For quantum double products set $z=y$. For quantum products set
$y=z=0$ and keep only $t+b=k$.
\end{method}

\begin{proposition}\label{prop:qcorrect}
Method~\ref{method:qproduct} computes \eqref{eq:qproduct}.
\end{proposition}
\begin{proof}
Multiply \eqref{eq:qexpansion} by $\widetilde\Sch_u(x;y)$ and apply
\eqref{eq:qpieri} at each column, with $p'=k-b$.
The coefficients in \eqref{eq:qpieri} involve only $q,y,z$, so they
commute with the remaining factors.
The proof of Theorem~\ref{thm:correct} applies verbatim, and no
contribution is lost since $p'\leq p_w$.
\end{proof}

\subsection{States of the quantum triple computation}

Write $D(w)=r_w$ for the last descent.
Put
\[
 \rho=\max(r_u,\tilde h),\qquad
 s(d)=\max\{j\geq0:j(j+1)\leq d\},\qquad
 R_q=\min\bigl(M-1,\ \rho+\min(L,s(d))\bigr),
\]
and
\[
 U_q=\min\left\{\fall{M}{R_q},\binom{d+R_q}{R_q}\right\},
 \qquad M=n+L.
\]

\begin{lemma}[Quantum descent steps]\label{lem:qdescent}
Let $t_{ab}$ be a quantum cover of $w$ with $a\leq k<b$.
If it is a Bruhat cover then $D(wt_{ab})\leq\max(D(w),k)$.
If it is length-decreasing then $b\leq D(w)+1$ and
$D(wt_{ab})\leq\max(D(w),b)$.
In particular, if $D(w)\geq k$, a single quantum cover raises $D$ by at
most one.
\end{lemma}
\begin{proof}
Put $w'=wt_{ab}$ and $D=D(w)$.
The entries of $w$ after position $D$ increase.
Only positions $a$ and $b$ change, so a descent of $w'$ after
$\max(D,k)$ can only occur at $b-1$ or $b$ with $b>D$.

For a Bruhat cover, $w'(b)=w(a)<w(b)<w(b+1)$, so there is no descent at
$b$.
If $b-1=a$ there is no descent at $b-1$ because $w'(a)=w(b)>w(a)$.
Otherwise $w(b-1)<w(b)$, and the cover condition forces
$w(b-1)<w(a)=w'(b)$.

For a length-decreasing cover, suppose $b\geq D+2$.
If $a<b-1$ then $w(b-1)$ lies between $w(b)$ and $w(a)$, so
$w(b-1)>w(b)$; but $b-1>D$.
If $a=b-1$ then $a>D$ and $w(a)<w(b)$.
Either case is a contradiction, so $b\leq D+1$.
The entries after position $b$ are unchanged and increasing,
so $D(w')\leq\max(D,b)$.
\end{proof}

\begin{lemma}[Quantum state bounds]\label{lem:qstates}
Every permutation produced by a $\mathtt{qpieri}$ call in
Method~\ref{method:qproduct} lies in $S_M$, has length at most $d$,
and has no descent after $R_q$.
There are at most $U_q$ such permutations.
Every monomial $q^\delta$ that occurs involves only
$q_1,\ldots,q_{M-1}$ and has $|\delta|\leq d/2$.
\end{lemma}
\begin{proof}
\emph{Ambient size.}
Let $w\in S_N$ and $k<N$.
For $c>N$, $w(c)=c$.
Hence a first transposition with $b\geq N+2$ is not a quantum cover:
the value $N+1$ at position $N+1$ lies strictly between $w(a)$ and
$w(b)=b$, which is forbidden for a Bruhat cover, and $w(a)<w(b)$
rules out a length-decreasing cover.
Since the values $b$ weakly decrease, the call stays in $S_{N+1}$.
At column $i$ the height is $k\leq\tilde h\leq n-1$, so induction
gives $S_{n+i}\subseteq S_M$.
In particular, every quantum factor $q_a\cdots q_{b-1}$ has
$b\leq M$.

\emph{Length and $q$-degree.}
Classical covers add one to $\ell$; a length-decreasing cover $t_{ab}$
changes $\ell$ by $1-2(b-a)$ and multiplies by a monomial of degree
$b-a$.
Hence, along any chain of covers starting at $u$,
\begin{equation}\label{eq:qbalance}
 \ell(w)+2|\delta|=\ell(u)+T,
\end{equation}
where $T$ is the total number of transpositions.
Consider a sequence of pairs reaching a key $(w,a)$ at boundary $i$, with
chain lengths $t_j$ and label edge lengths $b_j$.
Then $t_j\leq\tilde\theta_j-b_j$, and the labels are length-additive,
so
\[
 T=\sum_{j\leq i}t_j\leq|\tilde\theta_{\leq i}|-\ell(a),\qquad
 |\tilde\theta_{\leq i}|=\sum_{j\leq i}\tilde\theta_j.
\]
The label $a$ lies on a complete label path to $\tilde a_*$, whose
later edges have total length at most $\sum_{j>i}\tilde\theta_j$.
Lemma~\ref{lem:strict} gives
$\ell(\tilde a_*)=|\tilde\theta|-\ell(v)$, so
$\ell(a)\geq|\tilde\theta_{\leq i}|-\ell(v)$ and $T\leq\ell(v)$.
A node produced from $w$ in column $i+1$ has $t\leq p_w=k-b$ for an
edge $a\to a'$ leaving a label paired with $w$.
Concatenating with a sequence reaching $(w,a)$ gives
$T+t\leq|\tilde\theta_{\leq i+1}|-\ell(a')\leq\ell(v)$, by the same
argument at boundary $i+1$.
By \eqref{eq:qbalance} that node has length at most $d$, and the
monomial along this sequence has degree at most $d/2$.
All monomials attached to the node come from such sequences.

\emph{Descents.}
Fix such a sequence ending at a node $w$.
All heights satisfy $k\leq\tilde h\leq\rho$, and $D(u)\leq\rho$.
By Lemma~\ref{lem:qdescent}, a cover applied to a permutation with
$D\leq\rho+j-1$ produces one with $D\leq\rho+j$, and the bound
$D\leq\rho+j$ is reached from $D=\rho+j-1$, for $j\geq1$, only by a
length-decreasing cover with $b=\rho+j$.
Since $a\leq k\leq\rho$, that cover has $q$-degree $b-a\geq j$.
If $D(w)=\rho+e$ with $e>0$, the sequence therefore contains covers
of $q$-degrees at least $1,2,\ldots,e$, so $e(e+1)/2\leq|\delta|\leq d/2$.
Furthermore, suppose a call first raises $\max(D,\rho)$ to $\rho+j$.
This happens through a cover with $b=\rho+j$.
All later covers in that call have $b\leq\rho+j$, and by
Lemma~\ref{lem:qdescent} they keep $D\leq\rho+j$.
Hence each call raises $\max(D,\rho)$ by at most one.
A sequence through column $i$ makes $i\leq L$ calls.
This proves $D(w)\leq R_q$; the bound $D(w)\leq M-1$ is automatic.

The count $U_q$ follows as in Lemma~\ref{lem:states}.
\end{proof}

\begin{proof}[Proof of Theorem~\ref{thm:quantum}]
Correctness is Proposition~\ref{prop:qcorrect}.
The cost is counted exactly as in the proof of Theorem~\ref{thm:main},
with $U,V,h,\theta$ replaced by $U_q,\widetilde V,\tilde h,\tilde\theta$,
using Lemmas~\ref{lem:strict} and \ref{lem:qstates} and the strict
label bound above.
The table has at most $U_q\widetilde V$ keys because
monomials in $q$ are part of the coefficient.
The ambient and degree assertions are in Lemma~\ref{lem:qstates}.
For the coarse bound, use $L\leq\tilde h\leq n-1$, $M\leq2n-1$,
$U_q\leq(2n)^{R_q}$, $\widetilde V\leq n^{\tilde h}$,
$Q_{\tilde h}(n)\leq2n^{\tilde h}$, and $R_q\leq\rho+L\leq2\rho$.
The final claims follow from \eqref{eq:strictheight}.
\end{proof}

\subsection{Remarks on the bounds in the quantum setting}

\paragraph{Dependence on \texorpdfstring{$L$}{L}}
Classically $h=r_v$, so fixing code lengths suffices.
In Method~\ref{method:qproduct} the factor $\tilde h$ also requires
$L$ to be bounded.
This is not a weakness of the state estimate alone: the example
$\Sch_v=x_1^m$ has $r_v=1$ and $\tilde h=m$.
The difference comes from the distinct-height condition in
Theorem~\ref{thm:qexpansion}.
Both $L$ and $\tilde h$ are computed from the input code of $v$, and
$L\leq\min(\ell(v),n-1)$.
For quantum Schubert multiplication one may interchange $u$ and $v$.

\paragraph{Pessimistic worst-case parameters}
For $u=v=(m+1)12\cdots m$, $2\leq m\leq8$, the measured numbers of
distinct $u$-side states per column were $m^2+m+1$, and the label sets
had $m$ elements. The corresponding bounds $U_q$ and $\widetilde V$ grow
superpolynomially.
We do not know whether polynomial bounds in $n$, for fixed code
lengths, hold for the strict label graphs.

\paragraph{Medium-shape implementation}
The package's fast quantum and quantum double kernels use the medium shape
$\bar\theta$ of \cite[Definition~4.1]{samuelquantum}, which allows
certain pairs of equal heights and has $\bar\theta_1\leq\tilde h$
\cite[Proposition~4.3]{samuelquantum}.
Pairs are processed by the doubled Pieri rule
\cite[Methods~4.14--4.15]{samuelquantum}.
The termwise-quantization argument also applies to these shapes.
Write $X_{k+1}=X_k\cup\{x_{k+1}\}$.
Splitting $e_c(X_{k+1})=e_c(X_k)+x_{k+1}e_{c-1}(X_k)$ and reindexing
proves
\[
 e_a(X_k)e_b(X_k)
 =\sum_{c\leq b}e_c(X_{k+1})e_{a+b-c}(X_k)
  -\sum_{c>a}e_c(X_{k+1})e_{a+b-c}(X_k).
\]
In a medium shape the column before a pair of height $k$ has height
at least $k+2$, so the heights $k,k+1$ are not used elsewhere.
Thus $\mathcal Q$ maps each filtered medium summand to the product of
its single quantized factors and its doubled factors
\cite[(4.1)]{samuelquantum}.
This gives a filtered version of \cite[Theorem~4.7]{samuelquantum}.
A doubled step must, however, retain the intermediate permutation of
its two chains, so the same accounting gives a bound with $U_q^3$ in
place of $U_q^2$.
We do not analyze the implementation's additional degree cap.

\paragraph{Prior work}
Quantum Schubert polynomials \cite{fomin1997quantum}, the quantum Pieri
rule \cite{postnikov1999quantum}, and quantum double Schubert polynomials
with their equivariant quantum interpretation
\cite{mihalcea2006equivariant,lam2014quantum} give algorithms
for these coefficients.
We are not aware of a previous explicit bound for computing an entire
product in this quantum setting on complete flag varieties in terms of
input parameters.
Since the $q^0$ part of \eqref{eq:qproduct} is \eqref{eq:product},
the hardness remarks in Section~\ref{sec:prior} apply to it.

\section{Implementation and computational checks}\label{sec:implementation}

The package's ordinary and double kernels use the same main data flow,
as Proposition~\ref{prop:uniform} predicts:
compute the dominant code of $v^{-1}$, construct backward label graphs,
enumerate Pieri endpoints once per first-component state, and
accumulate matched contributions.
Methods~\ref{method:labels} and \ref{method:product} follow the classical
algorithm of \cite[Methods~4.9--4.11]{samuelquantum}.
In particular, the package's \texttt{kdown\_perms} restriction to positions
still holding their original value is the distinct-position condition
justified by Lemma~\ref{lem:cancellation}.
The Cauchy filter in Lemma~\ref{lem:filter} separately explains which
column endpoints may be passed to the next layer.

The remaining distinction concerns arithmetic evaluation.
The package's elementary-weight routine uses an unmemoized
divide-and-conquer expression construction.
The circuit implementation analyzed in Theorem~\ref{thm:main} uses
the shared recurrence
\eqref{eq:weightrec}.
For a literal recursive implementation, a conservative upper bound
is obtained by multiplying \eqref{eq:main} by
$(h+1)^{O(\log(h+1))}$.
This leaves the general $2^{O(n\log n)}$ bound and fixed-code-length
polynomiality unchanged, but is not the same operation count.

Arithmetic-circuit complexity is also not a wall-clock bound for
a computer algebra system that canonicalizes expressions, copies
subtrees, or performs automatic simplification.
Probabilistic removal of zero coefficients is not part of the
algorithm analyzed here.
Theorems~\ref{thm:main} and \ref{thm:quantum} apply to the explicitly
specified algorithms: the triple Method~\ref{method:product} and the
strict-shape quantum triple Method~\ref{method:qproduct}.
The former includes ordinary and double products; the latter includes
quantum and quantum double products.
They do not cover Grothendieck kernels, the independent transition--Monk
kernel, the medium-shape quantum and quantum double kernels
with their doubled Pieri steps (see Section~\ref{sec:quantum}),
merged-column implementations, or every released optimization.
The package's reference quantum and quantum double routines use the strict
shape and single-column steps, as Method~\ref{method:qproduct} does.

\subsection{Parabolic quantum and parabolic quantum double output}\label{sec:parabolic-output}

The parabolic quantum and parabolic quantum double rings of
Section~\ref{sec:dictionary} have distinct products and polynomial
representatives from the corresponding complete-flag rings.
The implementation computes their coefficients by first computing a
product in the corresponding complete-flag ring and then applying the
Peterson--Woodward comparison term by term. This coefficient comparison
is separate from the final finite-rank quotient.
The correctness of this reduction, including the quantum triple coefficients, is
proved in \cite[Section~9]{samuelquantum}. The implementation does not
construct a lift from a requested parabolic degree. Instead, it splits
each complete-flag coefficient into its $q$-monomials, tests whether each
exponent vector is a Peterson--Woodward lift, and, for an accepted
monomial, multiplies the indexing permutation by the two parabolic
longest elements, tests minimality, and reindexes the surviving quantum
variables.

This gives a direct runtime accounting for the projection. Let $N$ be
the ambient rank used by the projection, let $C$ be the number of
permutation--$q$-monomial pairs presented to it, and let $F$ be the work
of splitting the coefficient expressions by their $q$-monomials.
For each such pair, the exponent tests and reindexing take
$O(N^2)$ elementary operations. The current \texttt{check\_blocks}
implementation checks all subinterval sums by nested loops, and the
longest-element routine constructs the relevant elements by adjacent
swaps. Each therefore costs at most $O(N^3)$ operations. Permutation
multiplication, the minimality test, and the ambient-rank test cost at
most $O(N^2)$ in this bound. Consequently, the projection costs
\[
 O(F+N^3 C).
\]
Here $F$ is kept separate because \texttt{factor\_out\_q} recursively
decomposes symbolic expressions and convolves the resulting monomial
maps at multiplication nodes; its cost is not one operation per
surviving monomial.

In the path-counting model, each distinct permutation--$q$-monomial
pair in the complete-flag output has at least one path contribution, so
$C$ is bounded by the work of the complete-flag computation. Thus the
comparison adds at most a polynomial-in-$N$ factor to that path-count
bound, in addition to the explicit symbolic $q$-splitting work $F$.
This is a postprocessing bound for the implemented reduction, not a
faster multiplication algorithm: the complete-flag product is still
computed first, and the projection does not reduce its state space.

\subsection{Computational checks}

The specified update was compared with the
ordinary Python kernel on all $616$ pairs from $S_2,S_3,S_4$.
The double update was checked on all pairs from $S_3$ at three
integer specializations and on $40$ additional pairs from $S_4$,
both against the double Python kernel and using independent
divided-difference evaluation of the Schubert polynomials on both
sides of \eqref{eq:product}.
The dominant-completion and descent assertions were also checked
through $S_6$, and the boundary counts on $528$ weak intervals
through $S_5$.
The label graphs and their coefficient alphabets were compared with
the implementation on all $148$ nonidentity inputs from $S_2$ through $S_5$.
The distinct-position enumeration was compared with
\texttt{kdown\_perms}, and with the unrestricted weighted path sum,
on $2728$ column instances through $S_5$.
The elementary recurrence was checked in all $91$ degrees
$0\leq p\leq m\leq12$, including repeated coefficient parameters.
For Corollary~\ref{cor:grass}, the rectangular completion and the
Grassmannian property of every enumerated $u$-side state were checked
on $4233$ pairs of Grassmannian permutations through $S_7$.

For Section~\ref{sec:quantum}, Method~\ref{method:qproduct} was implemented
independently, with $\mathtt{qpieri}$ written from the definition of a
quantum cover and compared with \texttt{elem\_sym\_perms\_q} on every
call.
Its quantum specialization ($y=z=0$) agreed with both package quantum kernels,
the strict-shape reference and the medium-shape fast kernel, on all
$584$ pairs with $v\neq\id$ from $S_2,S_3,S_4$ and on $149$ random
pairs from $S_5$.
The quantum triple coefficients agreed with both quantum double kernels at
random integer specializations of $y,z,q$ on all $30$ such pairs from
$S_3$ and $60$ random pairs from $S_4$.
In the quantum double specialization ($z=y$), the computed coefficients satisfied
$c_{u,v}^{w,\delta}=c_{v,u}^{w,\delta}$ on $555$ pairs through $S_4$.
In all these runs the assertions of Lemma~\ref{lem:qstates} and the label
bound $|\tilde A_i|\leq\widetilde V$ were checked for every
enumerated state.
Lemma~\ref{lem:strict} was checked on all $5906$ nonidentity
permutations through $S_7$.

\section{Practical efficiency}\label{sec:practical}

The worst-case bounds do not predict the runtime or returned dictionary
size of an individual product. We therefore complement
them with measurements of the public ring interface. These experiments
use the default compiled kernels, including the medium-shape quantum
and quantum double kernels, rather than a literal implementation of the
arithmetic-circuit algorithms in Theorems~\ref{thm:main} and \ref{thm:quantum}.
The purpose is to identify costs that matter in use, not to infer
asymptotic optimality from a timing table.

\subsection{Workloads motivated by applications}

Two common tasks are the expansion of products in a Schubert basis
and the extraction of intersection numbers from repeated products.
The Sage \texttt{lrcalc} documentation gives ordinary Schur and quantum
Grassmannian examples \cite{sagelrcalcdocs}; the \texttt{Schubert2}
documentation computes the number of lines meeting four general lines
in projective three-space by integration on $\operatorname{Gr}(2,4)$
\cite{schubert2integral}. These motivate the following tests.

\paragraph{Intersection numbers}
For $\operatorname{Gr}(k,n)$, we multiply by the divisor class
$\sigma_1$ successively $k(n-k)$ times, imposing the finite-rank
quotient after each multiplication. The final coefficient is the
Pl\"ucker degree
\begin{equation}\label{eq:grassdegree}
 \int_{\operatorname{Gr}(k,n)}\sigma_1^{k(n-k)}
 =\frac{(k(n-k))!}
        {\displaystyle\prod_{i=1}^{k}\prod_{j=1}^{n-k}(n-i-j+1)}.
\end{equation}
We use $(k,n)=(2,4),(3,6),(4,8),(5,10)$, checking each result
against the rectangular hook-length formula \eqref{eq:grassdegree}.
The first is the four-lines problem, whose answer is $2$.
This tests repeated basis multiplication and quotient handling; it
computes an intersection number, not the individual geometric solutions.

\paragraph{Full Schur expansions}
To separate a growing output from growing code length, we compute
$h_A(x_1,x_2)^2$ for $A=4,8,12$ in the stable ordinary ring.
The input code is $(0,A)$, and the result has exactly $A+1$ basis
terms, as in Section~\ref{sec:complexity}.
No finite-rank quotient is imposed in these tests.

\paragraph{A published quantum example}
Buch's Quantum Calculator gives the example
\begin{equation}\label{eq:qcalcbenchmark}
 \begin{split}
 \sigma_{(2,1)}^{\star3}
 ={}&4\sigma_{(4,4,1)}+8\sigma_{(4,3,2)}
       +2\sigma_{(3,3,3)}\\
    &+5q\,\sigma_{(2)}+4q\,\sigma_{(1,1)}
 \qquad\text{in }QH^*(\operatorname{Gr}(3,7))
 \end{split}
\end{equation}
\cite{buchqcalc}.
We reproduce it with \texttt{QPSx(3,4)}, using two left-associated
multiplications and \texttt{kill\_ideal()} after each one.
For a partition $\lambda$ padded to $k$ parts, the permutation used
by the Python interface is
$\operatorname{uncode}(\lambda_k,\ldots,\lambda_1)$.
The full coefficient dictionary, including both quantum terms, is
checked against \eqref{eq:qcalcbenchmark}. The timing includes the
Peterson--Woodward projection, not merely the complete-flag kernel.

\paragraph{A controlled full-flag family}
To compare coefficient rings at identical inputs, we also use
\begin{equation}\label{eq:benchmarkfamily}
 \code(u_m)=(m,0,1),\qquad \code(v_m)=(0,m,1),
 \qquad m\in\{2,4,6,8\}.
\end{equation}
Here $n=m+2$, $r_u=r_v=3$, and $d=2m+2$.
For each pair we compute six distinct products: ordinary, double, triple,
quantum, quantum double, and quantum triple. The triple and quantum triple
products use independent coefficient alphabets.
All outputs are stable products: labels beyond the input $S_n$ are
retained. This is a designed scaling experiment, not a random sample
of research workloads.

\subsection{Measurement protocol and checks}

The experiments were run on November 10, 2026, using the non-editable
PyPI release \texttt{schubmult==5.3.0} in a CPython~3.14.10 virtual
environment, with SymEngine~0.15.0 and SymPy~1.14.0, on an Intel Core
i7-13850HX under Linux on WSL2.
Workers ran serially, each pinned to logical CPU~0, with OpenMP and
BLAS thread counts set to one. The local compiled extension was enabled;
its SHA-256 digest is recorded with the measurements.
This was a shared workstation, not an isolated benchmarking host.

Each repetition starts a fresh interpreter, imports the package from
the virtual environment, and constructs the operands before timing.
This avoids measuring hits in the ring's product cache as though they
were new multiplications. Deferred work triggered inside a timed call
is included. Only multiplication is timed: we record wall and process-CPU
times separately and report median wall times over five repetitions.
Multiplication excludes process startup, initial imports, printing,
and verification.
For the intersection-number and parabolic tests, the quotient and
projection steps described above are included.

We record the number $K$ of dictionary keys returned by the package,
without expanding the coefficients or counting their monomials.
Thus $K$ describes the returned data structure, not necessarily the
mathematical support: a stored coefficient expression may represent
the zero polynomial. No probabilistic zero pruning is enabled.

Multiplication and the subsequent untimed validation phase each have
an operating-system-enforced $30$-second limit, including while a
compiled kernel is running; each worker has a $4$-GiB address-space
limit. A censored observation records its phase and limit, and stops
further repetitions of that case. It is not replaced by a completed
timing or silently omitted.

All correctness checks run after the multiplication timer has stopped.
For every full-flag product, setting $y=z=q=0$ directly in the
unexpanded coefficients is checked against Buch's independent C
\texttt{schubmult} executable from \texttt{lrcalc}~1.2.
For $m=2$, all six products are additionally checked by exact polynomial
reconstruction from the output basis expansion; polynomial expansion
is used only for this small, untimed correctness check, not as a benchmark.
The application-based tests use \eqref{eq:grassdegree}, the two-row
Schur identity, and \eqref{eq:qcalcbenchmark} as exact checks.
Specialization checks alone do not certify every symbolic coefficient
for the larger inputs.

The companion files \texttt{schubmult\_benchmarks.py} and
\texttt{schubmult\_benchmarks.json} contain the executable protocol,
all input specifications, individual timings, output counts, check results, and
environment metadata. For example:
\begin{verbatim}
python3.14 -m venv .venv
.venv/bin/python -m pip install schubmult==5.3.0 \
  symengine==0.15.0 sympy==1.14.0
.venv/bin/python schubmult_benchmarks.py --installed-version 5.3.0 \
  --reference /usr/bin/schubmult --repeats 5 --timeout 30 \
  --output schubmult_benchmarks.json
\end{verbatim}
The reference executable is used only for checking; these are not
comparative timings against \texttt{lrcalc}, Sage, \texttt{Schubert2},
or Quantum Calculator. In particular, we do not compare the C
executable's process-level time with the Python ring interface's
post-import multiplication time: those measurements include different
startup costs and would not isolate the multiplication kernels.

\subsection{Measured performance and practical guidance}

All $32$ cases completed all five repetitions and the stated correctness
checks within the resource limits, giving $160$ multiplication measurements.

Table~\ref{tab:practical-applications} gives the application-based
measurements. The four intersection numbers are respectively
$2$, $42$, $24024$, and $701149020$, although at most $2$, $3$, $8$,
and $20$ basis keys were retained at any quotient boundary.
Thus a large numerical answer need not require a large intermediate
basis expansion. The quantum Grassmannian example also has a small
output and a short median runtime, with the full parabolic reduction
included; its five multiplication times ranged from $2.97$ to $3.96$
milliseconds. These are timings after initial imports, not command-line
startup times, and small differences between individual measurements
should not be interpreted as algorithmic effects.

\begin{table}[htbp]
\centering
\caption{Application-based workloads. Times include every multiplication
and quotient or projection step in the specified computation.
$K$ counts returned dictionary keys; the last
column gives the minimum and maximum of the five wall times.}
\label{tab:practical-applications}
\begin{tabular}{lrrl}
\hline
Computation & $K$ & Median (ms) & Range (ms)\\
\hline
$\operatorname{Gr}(2,4):\ \sigma_1^4$ & 1 & 0.209 & $[0.191,0.220]$\\
$\operatorname{Gr}(3,6):\ \sigma_1^9$ & 1 & 0.448 & $[0.430,0.628]$\\
$\operatorname{Gr}(4,8):\ \sigma_1^{16}$ & 1 & 1.45 & $[1.34,1.64]$\\
$\operatorname{Gr}(5,10):\ \sigma_1^{25}$ & 1 & 5.12 & $[4.89,5.36]$\\
$h_4(x_1,x_2)^2$ & 5 & 0.136 & $[0.132,0.159]$\\
$h_8(x_1,x_2)^2$ & 9 & 0.210 & $[0.191,0.223]$\\
$h_{12}(x_1,x_2)^2$ & 13 & 0.323 & $[0.310,0.360]$\\
$\operatorname{Gr}(3,7):\ \sigma_{(2,1)}^{\star3}$ & 5 & 3.06 & $[2.97,3.96]$\\
\hline
\end{tabular}
\end{table}

Table~\ref{tab:practical-flags} compares the six full-flag multiplication
tasks at identical inputs.
At $m=8$, the ordinary product returns $25$ keys in a median of
$0.341$ milliseconds, while the classical triple product returns
$142$ keys in $8.32$ milliseconds. Quantum, quantum double, and quantum
triple multiplication return $48$, $4433$, and $5935$ keys, with median
times of $13.1$, $110$, and $181$ milliseconds, respectively.
The shared complexity bound of Proposition~\ref{prop:uniform} is
therefore not a claim of equal wall-clock costs. A specialized run
can visit fewer states, and a coefficient-ring operation can cost much
more than an integer addition.

\begin{table}[htbp]
\centering
\small
\caption{The full-flag family \eqref{eq:benchmarkfamily}.
O, D, T denote ordinary, double, and triple products; Q, QD, QT
denote quantum, quantum double, and quantum triple products, respectively.
$K$ counts returned dictionary keys, not necessarily nonzero coefficients.
$t_{\rm mul}$ is the median multiplication wall time over five repetitions,
in milliseconds.}
\label{tab:practical-flags}
\begin{tabular}{rlrr}
\hline
$m$ & Mode & $K$ & $t_{\rm mul}$\\
\hline
2 & O & 6 & 0.121\\
2 & D & 15 & 0.222\\
2 & T & 17 & 0.247\\
2 & Q & 11 & 0.211\\
2 & QD & 21 & 0.296\\
2 & QT & 23 & 0.302\\
\hline
4 & O & 13 & 0.164\\
4 & D & 46 & 0.777\\
4 & T & 48 & 0.888\\
4 & Q & 24 & 0.876\\
4 & QD & 116 & 2.55\\
4 & QT & 143 & 3.40\\
\hline
6 & O & 19 & 0.257\\
6 & D & 87 & 2.16\\
6 & T & 89 & 2.85\\
6 & Q & 36 & 3.84\\
6 & QD & 713 & 15.9\\
6 & QT & 915 & 26.4\\
\hline
8 & O & 25 & 0.341\\
8 & D & 140 & 5.53\\
8 & T & 142 & 8.32\\
8 & Q & 48 & 13.1\\
8 & QD & 4433 & 110\\
8 & QT & 5935 & 181\\
\hline
\end{tabular}
\end{table}

These measurements suggest concrete choices for users.
When only ordinary coefficients are needed, use the ordinary ring
rather than compute a symbolic product and specialize afterwards.
Likewise, use \texttt{QSx} for quantum products and \texttt{QDSx}
when equivariant quantum coefficients are required.
For repeated products in a fixed Grassmannian, impose the intended
quotient after each step rather than retain irrelevant stable terms.
Distinguish a raw dictionary length from an exact support count.
Finally, reuse already computed products when appropriate, while
keeping cache-hit timings separate from first-computation timings.

The experiments demonstrate useful small and moderate computations,
not a universal practical advantage or an empirical proof of a
complexity class. Specialized Grassmannian routines remain relevant,
and looking up one coefficient after constructing a full product
does not turn \alg\ into a single-coefficient algorithm.
The tables are best read together with the output model and the
input-parameter bounds, rather than as substitutes for either.

\section{Further questions}

Several refinements of the input-parameter bounds remain.
One is to replace $Q_h(n)$ by a smaller count exploiting the
Bruhat-cover conditions, beyond the distinct-position restriction
already justified by \cite{samuelquantum}.
For Grassmannian inputs this is the only obstruction to a $2^{O(n)}$
bound in Corollary~\ref{cor:grass}.
Another is to exploit the degree cap and the actual dominant column
heights simultaneously, instead of replacing all column heights by $h$.

For quantum, quantum double, and quantum triple products, the main
question is whether the dependence on $L$ is necessary.
The example $x_1^m$ in Section~\ref{sec:quantum} has quadratically many
$u$-side states and linearly many labels, while
Theorem~\ref{thm:quantum} allows exponentially many.
A bound polynomial in $n$ for fixed $r_u,r_v$ would require either a
quantization argument tolerating repeated column heights, such as the
doubled factors of the medium shape, or a direct count of strict label
graphs.
A sharper bound on the descent excess in Lemma~\ref{lem:qstates} would
also help; in random tests through $S_9$ it never exceeded $2$.

The fixed-code-length bound also suggests a practical division of labor.
Explicit multiplication can generate whole tables for moderate numerical
code entries, while the bounded-range lattice-point method of
Pak and Slonim is designed for compressed single-coefficient inputs.
No asymptotic crossover or universal practical superiority follows from
the two available upper bounds alone.

\appendix
\section{Using the conjectural Grothendieck interface}\label{app:grothendieck}

The package provides \texttt{Gx} and \texttt{DGx} for ordinary and double
Grothendieck polynomials, and \texttt{QDGx} for quantum double
Grothendieck polynomials. Multiplication in \texttt{QDGx} uses a
conjectural equivariant quantum $K$-Pieri rule. Its non-equivariant
specialization at $\beta=-1$ is the Pieri rule of Naito and
Sagaki~\cite{naito2024pieri}. The Grothendieck kernels are outside the
scope of the complexity bounds proved in this paper.

\subsection{Python interface and conventions}

The constructor \texttt{QDGx(w, beta=-1)} accepts a permutation in
one-line notation and uses the coefficient alphabet $y$.
Its second argument selects another alphabet:
\texttt{QDGx(w, z, beta=-1)} uses $z$, while
\texttt{QDGx(w, "0", beta=-1)} gives the non-equivariant version.
Omitting \texttt{beta} leaves the deformation parameter symbolic;
$\beta=-1$ is the usual $K$-theory convention.
The double Grothendieck convention uses
$a\oplus b=a+b+\beta ab$. Consequently, at $\beta=0$ the quantum double
basis specializes to $\Sch_w^q(x;-y)$, with the opposite coefficient
alphabet from the Schubert convention used in the main text.

Elements support linear combinations, multiplication, and coefficient
access by \texttt{Permutation} keys. The method \texttt{simplify()}
cancels rational coefficients while retaining the quantum
Grothendieck basis; \texttt{expand()} instead returns the explicit
$x$-polynomial using Lenart--Maeno quantization~\cite{lenart2006quantum}.
For example:
\begin{verbatim}
>>> from schubmult import QDGx, Permutation
>>> from schubmult.abc import z
>>> g = QDGx([2, 1], "0", beta=-1)
>>> g.expand()
q_1 + x_1 - x_1*q_1
>>> p = (g * g).simplify()
>>> for w in sorted(p):
...     print(w, p[w])
() q_1
(1, 3, 2) -q_1
(3, 1, 2) 1
\end{verbatim}
Here \texttt{()} indexes the identity. These are stable polynomial
products: input permutations in $S_n$ do not impose a finite-flag
quotient, so output labels may lie outside $S_n$.
Conversion of arbitrary $x$-polynomials to the quantum Grothendieck
basis and coproducts are not implemented by \texttt{QDGx}.

Independent coefficient alphabets give quantum triple products.
Both factors must use the same value of $\beta$; the result is in the
left factor's coefficient alphabet. We print the coefficient dictionary
to keep the rational expressions readable:
\begin{verbatim}
>>> g = QDGx([2, 1], beta=-1)
>>> h = QDGx([2, 1], z, beta=-1)
>>> p = (g * h).simplify()
>>> for w in sorted(p):
...     print(w, p[w])
() -q_1*(-1 + z_1)/(1 - y_2)
(1, 3, 2) q_1*(-1 + z_1)/(1 - y_2)
(2, 1) (-y_2 + z_1)/(1 - y_2)
(3, 1, 2) (1 - z_1)/(1 - y_2)
\end{verbatim}
Specializing $z=y$ gives the equivariant product, and setting both
coefficient alphabets to zero gives the non-equivariant product.

\subsection{Multiplicative coefficient variables}

At $\beta=-1$, \texttt{ChangeOfVariables} can replace the default
coefficient coordinates by $1-y_i^{-1}$, so that the variables $y_i$
are multiplicative torus characters:
\begin{verbatim}
>>> from schubmult import ChangeOfVariables, QuantumDoubleGrothendieckRing
>>> from schubmult.abc import x, y
>>> characters = ChangeOfVariables(z, y, lambda i, t: 1 - 1/t[i])
>>> R = QuantumDoubleGrothendieckRing(x, characters, beta=-1)
>>> r = (R([2, 1])**2).simplify()
>>> r[Permutation([])]
y_2*q_1/y_1
>>> r[Permutation([2, 1])]
(y_1 - y_2)/y_1
\end{verbatim}
In these coordinates the non-equivariant specialization is $y_i=1$,
rather than $y_i=0$. No change of quantum parameters is required.

\subsection{Finite flags and parabolic projection}

For parabolic quantum $K$-theory, the Python interface applies
\texttt{apply\_kato} to the complete-flag product dictionary.
At $\beta=-1$ this implements Kato's quotient
map~\cite{kato2019quantum}: it replaces a permutation by its minimal
right-coset representative and sets the quantum parameters belonging
to the parabolic subgroup to $1$.
Unlike the command-line block sizes, the Python argument lists the
\emph{1-indexed simple reflections} generating that subgroup.
The argument \texttt{n} specifies the ambient $S_n$; terms indexed
outside $S_n$ are discarded, and the surviving quantum parameters
are reindexed consecutively. An empty generator list with an explicit
\texttt{n} gives the finite complete-flag restriction.

For blocks $(2,2)$, the generators are $[1,3]$ and $n=4$.
The complete-flag parameter $q_2$ becomes the parabolic parameter $q_1$.
Here is the mixed-alphabet analogue of the Grassmannian example in
Section~\ref{sec:examples}:
\begin{verbatim}
>>> from schubmult.mult.groth_quantum_double import apply_kato
>>> a = QDGx([3, 4, 1, 2], beta=-1)
>>> b = QDGx([1, 3, 2], z, beta=-1)
>>> coefficients = apply_kato(a * b, [1, 3], n=4, beta=-1)
>>> parabolic = a.ring.from_dict(coefficients).simplify()
>>> for w in sorted(parabolic):
...     print(w, parabolic[w])
(1, 3, 2) q_1*(-1 + z_1)*(-1 + z_2)/((1 - y_3)*(1 - y_4))
(3, 4, 1, 2) (-y_3 - y_4 + z_1 + z_2 + y_4*y_3 - z_1*z_2)/((1 - y_3)*(1 - y_4))
\end{verbatim}
The inputs are minimal coset representatives. Here \texttt{from\_dict}
is used only for simplification and display, not to construct a
parabolic ring. Use the same $\beta$ for multiplication and projection.
Since individual terms may contain inverse powers of $\beta$,
a $\beta=0$ specialization must follow symbolic projection and
coefficient cancellation.

\subsection{Command-line usage}

The four programs are:
\begin{center}
\begin{tabular}{ll}
\texttt{grothmult\_py} & Ordinary Grothendieck products\\
\texttt{grothmult\_double} & Double and mixed products\\
\texttt{grothmult\_q} & Non-equivariant quantum products\\
\texttt{grothmult\_q\_double} & Equivariant and mixed quantum products
\end{tabular}
\end{center}
All accept \texttt{--beta} to specialize the deformation parameter,
and \texttt{--code} to use Lehmer codes.
The double programs accept \texttt{--mixed-var} for independent
coefficient alphabets and \texttt{--simplify} to cancel coefficients.
The quantum examples above can be computed as follows:
\begin{verbatim}
$ grothmult_q 2 1 - 2 1 --beta=-1
()  q_1
(1, 3, 2)  -q_1
(3, 1, 2)  1
$ grothmult_q_double 2 1 - 2 1 --mixed-var --beta=-1 --simplify
       ()  (q_1 - z_1*q_1)/(1 - y_2)
(1, 3, 2)  (-q_1 + z_1*q_1)/(1 - y_2)
   (2, 1)  (-y_2 + z_1)/(1 - y_2)
(3, 1, 2)  (1 - z_1)/(1 - y_2)
\end{verbatim}
Parabolic inputs must be minimal coset representatives.
The quantum programs take \emph{block sizes} after
\texttt{--parabolic} and handle the delayed $\beta=0$
specialization automatically.

The option \texttt{--multiplicative-variables} selects $\beta=-1$
and the coordinates $1-y_i^{-1}$ for either double program, also using
$1-z_i^{-1}$ when \texttt{--mixed-var} is present.
Thus the parabolic triple product above, in multiplicative
coordinates, is:
\begin{verbatim}
$ grothmult_q_double 3 4 1 2 - 1 3 2 --parabolic 2 2 \
>     --mixed-var --multiplicative-variables --simplify
   (1, 3, 2)  y_4*y_3*q_1/(z_1*z_2)
(3, 4, 1, 2)  (-y_4*y_3 + z_1*z_2)/(z_1*z_2)
\end{verbatim}
The quantum double command does not currently implement
\texttt{--display-positive}.

\section*{Funding}
The author was not supported by any funding provided for the purpose of
preparing this research article.

\section*{Declaration of competing interests}
The author declares no competing interests.

\section*{Declaration of generative AI and AI-assisted technologies in the manuscript preparation}
During the preparation of this manuscript, GitHub Copilot via the
Copilot SDK in VS Code was used to assist with language editing,
mathematical review, literature searches, and code and benchmark checks.
The author reviewed and edited the resulting content and takes full
responsibility for the content of this publication.

\bibliographystyle{acm}
\bibliography{schubmult2}
\end{document}